\documentclass[preprint]{elsarticle}

\usepackage{amssymb,url,tikz,amsmath,mathdots}
\usepackage[english]{babel}

\newtheorem{theorem}{Theorem}[section] 
\newtheorem{lemma}[theorem]{Lemma}     
\newtheorem{corollary}[theorem]{Corollary}
\newtheorem{proposition}[theorem]{Proposition}
\newproof{pf}{Proof}
\newcommand{\Rel}{\mathcal{R}}
\newcommand{\BM}{{\mathfrak{A}}}
\newcommand{\Pone}{\widetilde{P}}
\newcommand{\Ptwo}{\widehat{P}}
\newcommand{\Mone}{M\llap{$\widetilde{\phantom{D}}$}}
\newcommand{\Mtwo}{M\llap{$\widehat{\phantom{D}}$}}
\newcommand{\Done}{\widetilde{D}}
\newcommand{\Dtwo}{\widehat{D}}
\newcommand{\irel}[1]{\stackrel{#1}\sim}
\newcommand{\V}{\mathcal{V}}
\newcommand{\hatGamma}{\widehat{\Gamma}}
\newcommand{\hatV}{\widehat{\V}}
\newcommand{\hatK}{\widehat{K}}
\newcommand{\Fq}{{\mathbb F}_q}
\newcommand\diag{\mathop{\rm diag}\nolimits}
\newcommand\diam{\mathop{\rm diam}\nolimits}
\newcommand{\Hom}{\mathop{\rm Hom}\nolimits}
\newcommand{\Aut}{\mathop{\rm Aut}\nolimits}
\newcommand{\Sym}{\mathop{\rm Sym}\nolimits}
\newcommand{\gauss}[2]{\genfrac[]{0 pt}1{#1}{#2}}
\newcommand{\Z}{{\mathbb Z}}

\begin{document}

\begin{frontmatter}
\title{Double Covers of Symplectic Dual Polar Graphs}

\author[uw1]{G.\ Eric Moorhouse\corref{cor1}}
\ead{moorhous@uwyo.edu}

\author[uw1]{Jason Williford}
\ead{jwillif1@uwyo.edu}

\cortext[cor1]{Corresponding author}
\address[uw1]{Department of Mathematics, University of Wyoming, Laramie WY 82071 USA}

\begin{abstract}
Let $\Gamma=\Gamma(2n,q)$ be the dual polar graph of type $Sp(2n,q)$. Underlying this graph is a $2n$-dimensional vector space $V$ over a field $\Fq$ of odd order~$q$,
together with a symplectic (i.e.\ nondegenerate alternating bilinear) form $B:V\times V\to\Fq$. The vertex set of $\Gamma$ is the set $\V$ of all $n$-dimensional totally isotropic subspaces of~$V$.
If $q\equiv1$ mod~4, we obtain from $\Gamma$ a nontrivial two-graph $\Delta=\Delta(2n,q)$ on~$\V$
invariant under $PSp(2n,q)$. This two-graph corresponds to a double cover
$\hatGamma\to\Gamma$ on which is naturally defined a $Q$-polynomial $(2n+1)$-class association scheme on $2|\hatV|$ vertices.
\end{abstract}
\end{frontmatter}

\section{Introduction}\label{intro}

Association schemes~\cite{BI,BCN} were first defined by Bose and Mesner~\cite{BM} in the context of the design of experiments.
Philippe Delsarte used association schemes to unify the study 
of coding theory and design theory in his thesis~\cite{PD}, where he derived his well-known linear programming bound which has since found many applications 
in combinatorics. There he identified two types of association schemes which were of particular interest: the so-called 
$P$-polynomial and $Q$-polynomial schemes. 
Schemes which are $P$-polynomial are precisely those arising from distance-regular graphs, and are well studied.
In particular, much effort has gone into the classification of distance-transitive graphs, the $P$-polynomial schemes 
which are the orbitals of a permutation group; and it is likely that all such examples are
known. Also well-studied are the schemes which are both $Q$-polynomial and
$P$-polynomial. A well-known conjecture~\cite[p.312]{BI} of Bannai and Ito is the following: for sufficiently
large $d$, a primitive scheme is $P$-polynomial if and only if it is $Q$-polynomial.

Classification efforts for $Q$-polynomial schemes are far less advanced than in the $P$-polynomial case; in particular it is likely that more
examples from permutation groups are yet to be found. The $Q$-polynomial property has no known combinatorial characterization, making their study more difficult. 
However, the list of known examples (see~\cite{MMW,PW,VDMM}) indicates that these objects have interesting structure from the viewpoint
of designs, lattices, coding theory and finite geometry.

In this paper, we give a new family of imprimitive $Q$-polynomial schemes with an
unbounded number of classes. These schemes are formed by the orbitals of a group, giving a double cover of the scheme arising from the symplectic dual polar space graph. 
We note that only one other family of imprimitive $Q$-polynomial schemes with an
unbounded number of classes is known that is not $P$-polynomial, namely the bipartite doubles of the Hermitian dual polar space graphs, which are $Q$-bipartite and $Q$-antipodal.
The schemes in this paper are $Q$-bipartite, and have two $Q$-polynomial orderings.
Except when the field order $q$ is a square, the splitting field of these schemes is also
irrational.  We note that this is the only known family of $Q$-polynomial schemes with
unbounded number of classes and an irrational splitting field.
In the last section we give open parameters for hypothetical primitive $Q$-polynomial subschemes of this family. 

Our paper is organized as follows: Background material on Gaussian coefficients, two-graphs
and double covers of graphs, are covered in Sections \ref{gaussian}--\ref{graphs}. In
Section~\ref{dual_polar} we recall the standard construction of the symplectic dual polar
graph $\Gamma=\Gamma(2n,q)$. There we also introduce the Maslov index, which we use in Section~\ref{double_cover} to construct the double cover
$\hatGamma\to\Gamma$ when $q\equiv1$ mod~4. In Section~\ref{scheme} we construct a $(2n+1)$-class association scheme
$\mathcal{S}=\mathcal{S}_{n,q}$ from $\hatGamma$; and in Section~\ref{Q-poly} we show that $\mathcal{S}$ is
$Q$-polynomial. The $P$-matrix of the scheme is constructed in Section~\ref{P-matrix}.
A particularly tantalizing open problem is the question whether $\mathcal{S}$ is in general
the extended $Q$-bipartite double of a primitive $Q$-polynomial scheme; see Section~\ref{subscheme}.

\section{Gaussian coefficients}\label{gaussian}

For all integers $n,k$ we define the {\em Gaussian coefficient\/}
\[\gauss{n}k=\gauss{n}k_q=\begin{cases}
\frac{(q^n-1)(q^{n-1}-1)\cdots(q^{n-k+1}-1)}{(q^k-1)(q^{k-1}-1)\cdots(q-1)},
&\mbox{if\ }k\geqslant0;\\
0,&\mbox{if\ }k<0.
\end{cases}\]
In particular for $k=0$ the empty product gives $\gauss{n}0=1$. In later sections, $q$ will be a fixed prime power; but here we may regard $q$ as an indeterminate, so that
for $n\geqslant0$, after cancelling factors we find $\gauss{n}k\in\Z[q]$; and specializing to
$q=1$ gives the ordinary binomial coefficients $\gauss{n}k_1=\binom{n}k$. For general
$n\in\Z$ we instead obtain a Laurent polynomial in $q$ with integer coefficients, i.e.\ 
$\gauss{n}k\in\Z[q,q^{-1}]$, as follows from conclusion (ii) of the following.

\begin{proposition}\label{prop2.1}
Let $n,k,\ell\in\Z$. The Gaussian coefficients satisfy
\begin{itemize}
\item[(i)]$\gauss{n}k=q^k\gauss{n-1}k+\gauss{n-1}{k-1}=\gauss{n-1}k+q^{n-k}\gauss{n-1}{k-1}$;
\item[(ii)]$\gauss{-n}k=(-q^{-n})^k\gauss{n+k-1}k$;
\item[(iii)]$\gauss{n}k\gauss{k}\ell=\gauss{n}\ell\gauss{n-\ell}{k-\ell}$;
\item[(iv)]$\gauss{n}k=\gauss{n}{n-k}$ whenever $0\leqslant k\leqslant n$.\qed
\end{itemize}
\end{proposition}
Most of the conclusions of Proposition~\ref{prop2.1} are found in standard references such as~\cite{A}. However, our definition of $\gauss{n}k$ differs from the
standard definition found in most sources, which either leave $\gauss{n}k$ undefined for $n<0$, or define it to be zero in that case.
Our extension to all $n\in\Z$ means that the recurrence formulas~(i) hold for all integers $n,k$, unlike the `standard definition' which fails for $n=k=0$. Property~(i) plays a role in
our later algebraic proofs using generating functions. In further defense of our definition,
we observe that it has become standard to extend the definition of binomial coefficients
$\binom{n}k$ so that $\binom{-n}k=(-1)^k\binom{n+k-1}k$ (see e.g.~\cite[p.12]{A});
and (ii) naturally generalizes this to Gaussian coefficients. We further note that (iii) holds for all $n,k\in\Z$
whether one takes the standard definition of $\gauss{n}k$ or ours. The one advantage of
the standard definition is that it renders superfluous the extra restriction $0\leqslant k\leqslant n$
in the symmetry condition~(iv). The interpretation of $\gauss{n}k$ as the number of $k$-subspaces
of an $n$-space over $\Fq$ is valid for all $n\geqslant0$.

In Section~\ref{P-matrix} we will make use of the well-known generating polynomials
\[E_m(t)=\prod_{i=0}^{m-1}(1+q^i t)=\sum_{\ell=0}^\infty q^{\binom{\ell}2}\gauss{m}\ell t^\ell\quad\mbox{for\ }m=0,1,2,\ldots;\]
note that in the latter sum, the terms for $\ell>m$ vanish, yielding $E_m(t)\in\Z[q,t]$ (or
after specializing to a fixed prime power $q$, we obtain $E_m(t)\in\Z[t]$).
Here we see the usual binomial coefficient $\binom{\ell}2=\frac12\ell(\ell-1)$. In Section~\ref{P-matrix} we will make use of the following obvious relations:

\begin{proposition}\label{prop2.2}
For all $m\geqslant0$, the generating function $E_m(t)$ satisfies
\begin{itemize}
\item[(i)]$E_m(-q t)=\frac{1- q^mt}{1-t} E_m(-t)$;
\item[(ii)]$E_m(q^2t)=\frac{1+q^{m+1}t}{1+qt} E_m(qt)$; and
\item[(iii)]$E_m(r^3 t)=\frac{1+rq^mt}{1+rt} E_m(rt)$ where $r=\sqrt{q}$.\qed
\end{itemize}
\end{proposition}

\section{Two-graphs and double covers of graphs}
\label{graphs}

Here we describe the most basic connections between two-graphs and double covers of graphs; see \cite{M,Se,BCN,Ta} for more details. Our notation is chosen to conform to that used in subsequent sections.

Let $\V$ be any set. Denote by $\binom{\V}k$ the collection of all $k$-subsets of $\V$
(i.e.\ subsets of cardinality $k$). A {\em two-graph on $\V$\/} is a subset
$\Delta\subseteq\binom{\V}3$ such that for
every 4-set $\{x,y,z,w\}\in\binom{\V}4$, an even number, i.e.\ 0, 2 or 4, of the triples
$\{x,y,z\}$, $\{x,y,w\}$, $\{x,z,w\}$, $\{y,z,w\}$ is in $\Delta$. If $\Delta$ is a
two-graph on $\V$, then the complementary set of triples
$\overline\Delta=\bigl\{\{x,y,z\}\in\binom\V3:\{x,y,z\}\notin\Delta\bigr\}$ is also a
two-graph, called the {\em complementary two-graph on $\V$\/}.

A {\em graph\/} on $\V$ is a subset $\Gamma\subseteq\binom\V2$. Elements of
$\Gamma$ are called {\em edges\/}.
The {\em complete graph on $\V$\/} is the graph $K_\V$ with full edge set $\binom{\V}2$.
In general the complementary set of pairs
$\overline\Gamma=\bigl\{\{x,y\}\in\binom\V2:\{x,y\}\notin\Gamma\bigr\}$ is the
{\em complementary graph on $\V$\/}.

Every graph on $\V$ may be identified with a signing of the edges of the complete graph $K_\V$, i.e.\ a function $\sigma:\binom{\V}2\to\{\pm1\}$. Under this correspondence, the graph corresponding to $\sigma$ has as its edge set $\sigma^{-1}(1)=\bigl\{\{x,y\}\in\binom{\V}2:\sigma(x,y)=1\bigr\}$. (Here we abbreviate
$\sigma(\{x,y\})=\sigma(x,y)$.)

Given $\Gamma$ and $\sigma$ as above (which amounts to two graphs which may be
entirely unrelated except for sharing the same vertex set $\V$), we construct a new graph
$\hatGamma=\hatGamma_\sigma$ with vertex set $\hatV=\V\times\{\pm1\}$ and adjacency relation defined by
\[(x,\varepsilon)\sim(y,\varepsilon')\;\iff\;x\sim y\mbox{\ and\ }\varepsilon\varepsilon'=\sigma(x,y).\]
(Note that $(x,1)\not\sim(x,-1)$ since $\Gamma$ has no loops.)
The map $(x,\varepsilon)\mapsto x$ is a {\em double covering map\/} $\theta:\hatGamma\to\Gamma$, also called a {\em double cover\/} or simply a {\em cover\/}; and the {\em fibers\/}
of this map are the pairs $\theta^{-1}(x)=\{(x,1),(x,-1)\}$ where $x\in\V$.
(By definition, a {\em covering map\/} of graphs is a graph homomorphism $\theta:\hatGamma\to\Gamma$ such that for any vertex $x\in\Gamma$, the preimage of the neighborhood graph $\Gamma_x$ is isomorphic to a disjoint union of copies of $\Gamma_x$; see e.g.~\cite{GR}. `Double' refers to the condition that the covering map is
2-to-1.) We also say that the vertices $(x,1)$ and $(x,-1)$ are {\em antipodal\/}
with respect to the covering map. (Note that antipodal vertices must be at
distance${}\geqslant2$; but we deviate from common custom by
{\em not requiring\/} pairs of antipodal vertices to be at maximal distance $\diam\hatGamma$.) We denote
by $\zeta$ the transposition interchanging antipodal vertices:
$(x,1)\stackrel\zeta\leftrightarrow(x,-1)$. Denote by
$\Aut_\zeta\hatGamma\leqslant\Aut\hatGamma$
the subgroup consisting of all automorphisms of the graph $\hatGamma$ which preserve the
antipodality relation. In general, $\Aut_\zeta\hatGamma$ is the centralizer of $\zeta$
in the full automorphism group $\Aut\hatGamma\leqslant\Sym\hatV$; but in our case
we obtain equality $\Aut_\zeta\hatGamma=\Aut\hatGamma$ (see
Lemma~\ref{lemma5.4}).
Similarly, two covers $\theta_i:\hatGamma_i\to\Gamma$ of the
same graph $\Gamma$ (for $i=1,2$) are {\em equivalent\/} or {\em isomorphic\/}
if there is a graph isomorphism $\rho:\hatGamma_1\to\hatGamma_2$ which preserves
antipodality, i.e.\ $\theta_1\circ\rho=\theta_2$.

Given $\sigma:\binom{\V}2\to\{\pm1\}$ as above, for every triple $\{x,y,z\}\in\binom{\V}3$ we may define
\[\sigma(x,y,z)=\sigma(x,y)\sigma(y,z)\sigma(z,x)\in\{\pm1\}.\]
A triple $\{x,y,z\}\in\binom{\V}3$ is called {\em coherent\/} or {\em non-coherent\/}
according as $\sigma(x,y,z)=1$ or $-1$.
The set of all coherent triples forms a two-graph on $\V$, denote by $\Delta_\sigma$; and
the set of non-coherent triples gives the complementary two-graph $\overline{\Delta}_\sigma$.

Two sign functions $\sigma_1,\sigma_2:\binom{\V}2\to\{\pm1\}$ (or the corresponding graphs $\sigma_1^{-1}(1)$, $\sigma_2^{-1}(1)$ on $\V$) are {\em switching-equivalent\/}
in the sense of Seidel~\cite{Se} if there exists a map $f:\V\to\{\pm1\}$ such that $\sigma_2(x,y)=f(x)f(y)\sigma_1(x,y)$ for all $\{x,y\}\in\binom{\V}2$. We have
$\Delta_{\sigma_1}=\Delta_{\sigma_2}$ iff $\sigma_1$ and $\sigma_2$ are
switching-equivalent. Assuming this holds,
then the corresponding covers $\hatGamma_{\sigma_1}$ and $\hatGamma_{\sigma_2}$
are isomorphic via $(x,\varepsilon)\mapsto(x,f(x)\varepsilon)$.

In the special case of the complete graph
$\Gamma=K_\V$, the following three notions are equivalent (see \cite[\S1.5]{BCN}): two-graphs on $\V$,
switching classes of graphs on $\V$, and isomorphism classes of double covers of the complete graph $K_\V$. For example given a double cover $\widehat{K_\V}\to K_\V$, the
corresponding two-graph is obtained as follows (see \cite[p.488]{Se}): Each triple $\{x,y,z\}$ of distinct vertices in
$\V$ induces a triangle $K_{\{x,y,z\}}\subseteq K_\V$; and such a triple is coherent iff its preimage in $\widehat{K_\V}$ induces a pair of triangles, rather than a 6-cycle, in $\widehat{K}_\V$.

An {\em automorphism of a two-graph $\Delta$\/} is a permutation of the underlying
point set $\V$ which preserves the set of coherent triples. We now relate $\Aut\Delta$
to the group $\Aut_\zeta\hatK\leqslant\Aut\hatK$ defined above for the associated double cover
$\hatK\to K$, where we abbreviate the complete graph $K_\V=K$. The following is easy to verify (or see \cite[\S2]{Ta}, where this isomorphism is denoted $\widehat{G}/Z\cong G$):

\begin{proposition}\label{prop3.1}
The group $\Aut_\zeta\hatK$ acts naturally on $\Delta$, inducing the full automorphism group of $\Delta$.
The kernel of this action is the central subgroup $\langle\zeta\rangle$ of order~2; thus $(\Aut_\zeta\hatK)/\langle\zeta\rangle\cong\Aut\Delta$.
\end{proposition}

\section{Dual polar graphs of type $Sp(2n,q)$, $q$ odd}
\label{dual_polar}

Fix a finite field $\Fq$ of odd prime power order $q$; an integer $n\geqslant1$; a $2n$-dimensional vector space $V$ over $\Fq$; and a symplectic (i.e.\ nondegenerate alternating) bilinear form $B:V\times V\to\Fq$. The {\em symplectic group\/} $Sp(2n,q)$ consists of
all {\em(linear) isometries\/} of $B$, i.e.
\[Sp(2n,q)=\{g\in GL(V):B(x^g,y^g)=B(x,y)\mbox{\ for all\ }x,y\in V\}.\]
The group of all {\em (linear) similarities\/} of $B$ is
\begin{align*}
GSp(2n,q)&=\{g\in GL(V):\mbox{for some nonzero\ }\mu\in\Fq\mbox{\ we have}\\ &\qquad B(x^g,y^g)=\mu B(x,y)\mbox{\ for all\ }x,y\in V\};
\end{align*}
some other notations for this group are $GSp_n(q)$ in~\cite{KL} or $CSp_n(q)$ in~\cite[p.31]{BHR}.
Replacing $GL(V)$ by $\Gamma L(V)\cong GL(V)\rtimes\Aut\Fq$, the group of all semilinear transformations of $V$, we obtain the group $\Sigma Sp(2n,q)$ of
all {\em semi-isometries\/}, and
the group $\Gamma Sp(2n,q)$ of all {\em semi-similarities\/} of $B$, given by
\begin{align*}
\Sigma Sp(2n,q)&=\{g\in\Gamma L(V):\mbox{for some\ }\tau\in\Aut\Fq\mbox{\ we have}\\
&\qquad\quad B(x^g,y^g)=B(x,y)^\tau\mbox{\ for all\ }x,y\in V\}\\
&\cong Sp(2n,q)\rtimes\Aut\Fq;\\
\Gamma Sp(2n,q)&=\{g\in\Gamma L(V):\mbox{for some nonzero\ }\mu\in\Fq\mbox{\ and\ }\tau\in\Aut\Fq\\
&\qquad\quad\mbox{we have\ }B(x^g,y^g)=\mu B(x,y)^\tau\mbox{\ for all\ }x,y\in V\}\\
&\cong GSp(2n,q)\rtimes\Aut\Fq.
\end{align*}
The projective versions of these groups are
\begin{align*}PSp(2n,q)&=Sp(2n,q)/\langle-I\rangle,\\
PGSp(2n,q)&=GSp(2n,q)/Z,\\
P\Sigma Sp(2n,q)&=\Sigma Sp(2n,q)/\langle-I\rangle,\\
P\Gamma Sp(2n,q)&=\Gamma Sp(2n,q)/Z
\end{align*}
where the central subgroup $Z$ of order $q-1$ consists of all scalar transformations $v\mapsto\lambda v$
for $0\neq\lambda\in\Fq$. We have
\[[P\Gamma Sp(2n,q):P\Sigma Sp(2n,q)]=[PGSp(2n,q):PSp(2n,q)]=2\]
where the nontrivial coset in both cases is represented by $h\in GSp(2n,q)$
satisfying $B(u^h,v^h)=\eta B(u,v)$ and $\eta\in\Fq$ is a nonsquare.

Our choice of notation for these groups, while not universal, is intended to conform
reasonably with \cite{Atlas,KL}. The group $P\Gamma Sp(2n,q)$, for example, is denoted $PC\Gamma Sp_n(q)$ in \cite[p.31]{BHR}. It arises (see Theorem~\ref{theorem4.1}) as the full automorphism
group of the associated dual polar graph, which we now describe.

Denote by $\V$ be the collection of all maximal totally
isotropic subspaces with respect to $B$, i.e.\ 
\[\V=\{X\leqslant V:X^\perp=X\}\]
where by definition $X^\perp=\{v\in V:B(x,v)=0\mbox{\ for all\ }x\in X\}$. Members of
$\V$ are often called {\em generators\/}, and every $X\in\V$ has dimension~$n$.
Denote by $\Gamma=\Gamma(2n,q)$ the graph on $\V$ where two vertices
$X,Y\in\V$ are adjacent iff $X\cap Y$ has codimension~1 in both $X$ and $Y$. More generally, the distance between $X$ and $Y$ in $\Gamma$ is
$d(X,Y)=k\in\{0,1,2,\ldots,n\}$ where the subspace $X\cap Y$ has codimension $k$ in both $X$ and~$Y$. Let $\Gamma_k$ denote the graph of the
distance-$k$ relation on $\V$; i.e.\ $\Gamma_k$ has vertex set $\V$ and two vertices $X,Y\in\V$ are adjacent in $\Gamma_k$ iff $d(X,Y)=k$.
The graph $\Gamma_1=\Gamma$ is called the {\em dual polar graph of type\/} Sp$(2n,q)$. It is {\em distance
regular\/}: given any two vertices $X,Y$ in $\Gamma$ at distance $k\in\{0,1,2,\ldots,n\}$, the vertex $Y$ has $q^{\binom{n-k}2}\gauss{n}k$ neighbors $Z$ in $\Gamma$, of which
\begin{align*}
a_k&=q^k-1\mbox{\ are at distance $k$ from\ }X,\\
b_k&=q^{k+1}\gauss{n-k}1\mbox{\ are at distance $k+1$ from $X$, and}\\
c_k&=\gauss{k}1\mbox{\ are at distance $k-1$ from\ }X;
\end{align*}
see \cite[\S9.4]{BCN}.
The edges of $\Gamma_1,\Gamma_2,\ldots,\Gamma_n$ partition the non-identical pairs
on $\V$, viewed as the edges of the complete graph $K_\V$; and together with the identity relation
$\Gamma_0=\{(X,X):X\in\V\}$ we obtain an $n$-class association scheme on $\V$ (see Section~\ref{scheme}). This scheme is $P$-polynomial since $\Gamma$ is distance regular; see~\cite{BCN}.
\begin{theorem}\label{theorem4.1}
For $n\geqslant2$, the full automorphism group of $\Gamma=\Gamma(2n,q)$ is the group $P\Gamma Sp(2n,q)$ acting naturally on the projective space of $V$.
\end{theorem}

\begin{pf}
See \cite[p.275]{BCN} (where this group is however denoted $P\Sigma p(2n,q)$).\qed
\end{pf}
Note that when $n=1$, the dual polar graph $\Gamma(2,q)$ is simply the complete graph
$K_{q+1}$, whose full automorphism group is the symmetric group of degree $q+1$.

For use in Section~\ref{double_cover} we record the following well-known fact. Although it follows easily from the axioms of polar geometry (or of near polygons), in the interest of self-containment we include a proof.

\begin{lemma}\label{lemma4.2}
The `diamond' graph (as shown) is not an induced subgraph of the dual polar graph $\Gamma$.
\end{lemma}
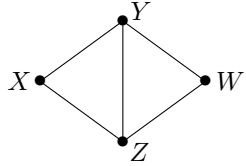
\[\begin{tikzpicture}[scale=1.1]
\node at (-1,0) {$\bullet$};
\node at (-1.25,0) {$X$};
\node at (0,.73) {$\bullet$};
\node at (.22,.85) {$Y$};
\node at (0,-.73) {$\bullet$};
\node at (.2,-.85) {$Z$};
\node at (1,0) {$\bullet$};
\node at (1.3,0) {$W$};
\draw (-1,0) -- (0,.73) -- (1,0) -- (0,-.73) -- (-1,0);
\draw (0,.73) -- (0,-.73);
\end{tikzpicture}\]

\begin{pf}If $X,Y,Z$ are mutually adjacent as shown, then $X\cap Y$ and $X\cap Z$ are distinct subspaces of codimension~1 in $X$, so $X=(X\cap Y)+(X\cap Z)$, whence $X\subseteq Y+Z$. Thus $X=X^\perp\supseteq Y^\perp\cap Z^\perp=Y\cap Z$. Similarly,
$W\supseteq Y\cap Z$. Now $X\cap W$ contains a subspace of dimension $n{-}1$, contradicting $d(X,W)\geqslant2$.\qed
\end{pf}

Now let $X$ be any $n$-dimensional vector space over $\Fq$. An $n$-linear form $f:X^n\to\Fq$ (i.e.\ linear in each argument whenever the other
$n-1$ arguments are fixed) is {\em alternating\/} if $f(x_1,x_2,\ldots,x_n)=0$ whenever two $x_i$'s coincide; equivalently,
$f(x_{1^\tau},x_{2^\tau},\ldots,x_{n^\tau})=-f(x_1,x_2,\ldots,x_n)$ for every odd permutation $\tau$ of the indices.
The space of all such alternating forms is one-dimensional, and is canonically identified with $(\bigwedge^{\mskip-2mu n}\mskip-2mu X)^*$,
the dual space of $\bigwedge^{\mskip-2mu n}\mskip-2mu X$. A {\em determinant function on $X$\/} is any nonzero alternating form $X^n\to\Fq$. Since
$\dim(\bigwedge^{\mskip-2mu n}\mskip-2mu X)^*=1$, a determinant function is determined up to nonzero scalar multiple.

Fix a choice of determinant function $\delta_X$ for each $X\in\V$. Although these choices are not canonical, one may proceed by arbitrarily choosing a basis
$\psi_1,\psi_2,\ldots,\psi_n$ for $X^*=\Hom(X,\Fq)$; then we obtain a determinant function on $X$ by defining
\[\delta_X(x_1,x_2,\ldots,x_n)=\det\bigl(\psi_i^*(x_j):1\leqslant i,j\leqslant n\bigr).\]

We need to define $\sigma(X,Y)\in\{\pm1\}$ for any pair $X\neq Y$ in $\V$. Let $k\in\{1,2,\ldots,n\}$ be the codimension of $X\cap Y$ in both
$X$ and $Y$. Choose bases $x_1,x_2,\ldots,x_n$ and $y_1,y_2,\ldots,y_n$ for $X$ and $Y$ respectively, such that $x_i=y_i$ (for $k<i\leqslant n$)
is a common basis for $X\cap Y$. (These bases depend on the choice of pair $(X,Y)$ and so are unrelated to any bases for $X$ and $Y$ used as a crutch
for constructing the corresponding determinant functions). Define
\[\sigma(X,Y)=\chi\bigl(\delta_X(x_1,x_2,\ldots,x_n)\delta_Y(y_1,y_2,\ldots,y_n)\det\bigl[B(x_i,y_j):1\leqslant i,j\leqslant k\bigr]\bigr)\]
where $\chi:\Fq^{\times}\to\{\pm1\}$ is the quadratic character:
$\chi(a)=1$ or $-1$ according as $a\in\Fq^{\times}$ is a square or a nonsquare.
This definition is implicit in~\cite{Th,KS}; and inspired by the literature, we refer to
$\sigma(X,Y)$ (or the ternary function $\sigma(X,Y,Z)$ defined below) as the
{\em Maslov index\/}. Note that $B$ induces a nondegenerate bilinear form on the
$2k$-space $(X+Y)/(X\cap Y)$, so that the $k\times k$
matrix $\bigl[B(x_i,y_j):1\leqslant i,j\leqslant k\bigr]$ is nonsingular.
\bigskip

\begin{proposition}\label{prop4.3}
Let $X,Y\in\V$ at distance $d(X,Y)=k\in\{0,1,2,\ldots,n\}$.
\begin{itemize}
\item[(i)]The value of $\sigma(X,Y)$ is independent of the choice of bases $x_i$ and $y_j$ as above.
\item[(ii)]Its dependence on the choice of determinant functions is expressed as follows:
Replacing $\delta_X$ by $c\delta_X$ has the effect of multiplying $\sigma(X,Y)$ by $\chi(c)$.
\item[(iii)]$\sigma(Y,X)=\chi(-1)^k\sigma(X,Y)=(-1)^{k(q-1)/2}\sigma(X,Y)$.
\item[(iv)]Let $g\in\Gamma Sp(2n,q)$, so that there exists a nonzero scalar $\mu_g\in\Fq$
and $\tau_g\in\Aut\Fq$ satisfying $B(x^g,y^g)=\mu_g B(x,y)^{\tau_g}$ for all $x,y\in Y$.
Then there exist nonzero scalars $\lambda_{g,U}\in\Fq$ for $U\in\V$, such that
\[\sigma(X^g,Y^g)=\chi\bigl(\mu_g^k\lambda_{g,X}\mskip1mu\lambda_{g,Y}\bigr)\sigma(X,Y).\]
\end{itemize}
\end{proposition}

\begin{pf}
Consider a change of basis on $X\cap Y$ specified by $x_i'=y_i'=\sum_{k<j\leqslant n}a_{ij}x_j$ where $A=\bigl(a_{ij}:k<i,j\leqslant n\bigr)$
is any invertible $(n-k)\times(n-k)$ matrix. Then
\[\delta_X(x_1,\ldots,x_k,x_{k+1}',\ldots,x_n')=(\det A)\delta_X(x_1,\ldots,x_n)\]
and $\delta_Y(y_1,\ldots,y_n)$ is multiplied by the same factor, $\det A$. The $(n-k)\times (n-k)$ matrix $[B(x_i,y_j):k<i,j\leqslant n]$ is unchanged, so the
value of $\sigma(X,Y)$ is multiplied by a net factor of $\chi\bigl((\det A)^2\bigr)=1$.

Next consider replacing $x_1,\ldots,x_k$ by $x_1',\ldots,x_k'$ where
\[x_i'\equiv\sum_{1\leqslant j\leqslant k}a_{ij}x_j\quad{\rm mod\ }(X\cap Y)\]
for $i=1,2,\ldots,k$ where $A$ is an invertible $k\times k$ matrix, and we leave the basis of $Y$ unchanged. Then
\[\begin{array}{c}
\delta_X(x_1',\ldots,x_k',x_{k+1},\ldots,x_n)=(\det A)\delta_X(x_1,\ldots,x_n);\\
\det\bigl[B(x_i',y_j):1\leqslant i,j\leqslant k\bigr]=(\det A)\det\bigl[B(x_i,y_j):1\leqslant i,j\leqslant k\bigr]
\end{array}\]
and the $\delta_Y$ factor is unchanged; so once again, the value of $\sigma(X,Y)$ is multiplied by $\chi\bigl((\det A)^2\bigr)=1$.
The same argument applies if $y_1,\ldots,y_k$ are replaced by $y_1',\ldots,y_k'$, and so~(i) follows. Conclusion~(ii) is clear.

Interchanging $X$ and $Y$ has the effect of interchanging the $\delta_X$ and $\delta_Y$ factors, and replacing
\begin{multline*}
\bigl[B(x_i,y_j):1\leqslant i,j\leqslant k\bigr]\mapsto\\
\bigl[B(y_i,x_j):1\leqslant i,j\leqslant k\bigr]=-\bigl[B(x_i,y_j):1\leqslant i,j\leqslant k\bigr].
\end{multline*}
The determinant of this matrix accrues a factor of $(-1)^k$, whence (iii) holds.

Let $g\in\Gamma Sp(2n,q)$. There exists a nonzero $\mu_g\in\Fq$ and
$\tau_g\in\Aut\Fq$ such that $(au+bv)^g=a^{\tau_g}u^g+b^{\tau_g}v^g$ and $B(u^g,v^g)=\mu_g B(u,v)^{\tau_g}$ for all
$a,b\in\Fq$ and $u,v\in V$. Now the map
\[X^n\to\Fq,\quad(x_1,x_2,\ldots,x_n)\mapsto\delta_{X^g}\bigl(x_1^g,x_2^g,\ldots,x_n^g\bigr)^{\tau_g^{-1}}\]
is a determinant function on $X$, so it is a scalar multiple of $\delta_X(x_1,x_2,\ldots,x_n)$. Hence there exists a nonzero scalar $\lambda_X=\lambda_{g,X}\in\Fq$ such that
\[\delta_{X^g}\bigl(x_1^g,x_2^g,\ldots,x_n^g\bigr)=\lambda_{g,X}\delta_X(x_1,x_2,\ldots,x_n)^{\tau_g}\]
for all $x_1,x_2,\ldots,x_n\in X$.

Now given $X,Y\in\V$ at distance $k$, fix bases $x_i,y_i$ as before; then
\begin{align*}
\sigma(X^g,Y^g)&=\chi\bigl(\delta_{X^g}(x_1^g,x_2^g,\ldots,x_n^g)\delta_{Y^g}(y_1^g,y_2^g,\ldots,y_n^g)\\
&{\mskip100mu}\times\det\bigl[B(x_i^g,y_j^g):1\leqslant i,j\leqslant k\bigr]\bigr)\\
&=\chi\bigl(\lambda_{g,X}\delta_X(x_1,x_2,\ldots,x_n)^{\tau_g}\lambda_{g,Y}\delta_Y(y_1,y_2,\ldots,y_n)^{\tau_g}\\
&{\mskip100mu}\times\det\bigl[\mu_g B(x_i,y_j)^{\tau_g}:1\leqslant i,j\leqslant k\bigr]\bigr)\\
&=\chi\bigl(\mu_g^k\lambda_{g,X}\lambda_{g,Y}\bigr)\chi\bigl(\delta_X(x_1,x_2,\ldots,x_n)^{\tau_g}\delta_Y(y_1,y_2,\ldots,y_n)^{\tau_g}\\
&{\mskip100mu}\times\det\bigl[B(x_i,y_j):1\leqslant i,j\leqslant k\bigr]^{\tau_g}\bigr)\\
&=\chi\bigl(\mu_g^k\lambda_{g,X}\mskip1mu\lambda_{g,Y}\bigr)\sigma(X,Y)
\end{align*}
since $\chi(a^\tau)=\chi(a)$. This proves (iv).\qed
\end{pf}

For each triple $(X,Y,Z)$ with distinct $X,Y,Z\in\V$, define
\[\sigma(X,Y,Z)=\sigma(X,Y)\sigma(Y,Z)\sigma(Z,X)\in\{\pm1\}.\]
A triple $(X,Y,Z)$ of distinct elements of $\V$ is {\em coherent\/} or {\em non-coherent\/}
according as $\sigma(X,Y,Z)=1$ or $-1$.

\begin{theorem} \label{theorem4.4}
Suppose $q\equiv1\,$ mod~$4$. Then the set of coherent
triples forms a two-graph $\Delta_\sigma$ on $\V$, invariant under $P\Sigma Sp(2n,q)$.
\end{theorem}

\begin{pf}
Let $X,Y,Z,W\in\V$ be distinct. Since $\chi(-1)=1$, $(X,Y,Z)$ is coherent iff any permutation of its members yields a coherent triple; so
the set of coherent triples may be regarded as a collection of unordered triples $\{X,Y,Z\}$. Since
\begin{multline*}
\sigma(X,Y,Z)\sigma(X,Y,W)\sigma(X,Z,W)\sigma(Y,Z,W)\\
=\sigma(X,Y)^2\sigma(X,Z)^2\cdots\sigma(Z,W)^2=1,
\end{multline*}
evenly many of the triples in $\{X,Y,Z,W\}$ are coherent. If $g\in\Gamma Sp(2n,q)$ with $B(x^g,y^g)=\mu_g B(x,y)^{\tau_g}$, then
\begin{align*}
\sigma(X^g,Y^g,Z^g)&=\chi(\mu_g^{d(X,Y)}\lambda_{g,X}\mskip1mu \lambda_{g,Y})\chi(\mu_g^{d(Y,Z)}\lambda_{g,Y}\mskip1mu\lambda_{g,Z})\\
&\phantom{{}={}}\times\chi(\mu_g^{d(Z,X)}\lambda_{g,Z}\mskip1mu \lambda_{g,X})\sigma(X,Y,Z)\\
\phantom{\sigma(X^g,Y^g,Z^g)}&=\chi(\mu_g)^{d(X,Y)+d(Y,Z)+d(Z,X)}\sigma(X,Y,Z).
\end{align*}
In particular when $g\in\Sigma Sp(2n,q)$, $\mu_g=1$ and $\sigma(X^g,Y^g,Z^g)=\sigma(X,Y,Z)$.\qed
\end{pf}

If $q\equiv3\,$ mod~$4$, or $g\in P\Gamma Sp(2n,q)$ with $g\notin P\Sigma Sp(2n,q)$,
the situation is a little trickier: various subsets of the coherent triples form either two-graphs or skew two-graphs in the sense of \cite{M},
invariant under $Sp(2n,q)$. We ignore this case here, and {\em henceforth assume that\/}
\[q\equiv1\,\hbox{\ mod\ }4.\]

We next show that {\em in a geodesic path, every triple of vertices is coherent.\/}

\begin{lemma} \label{lemma4.5}
Suppose $q\equiv1\,$ \hbox{mod $4$}. Let $X,Y,Z\in\V$ such that $d(X,Y)=j$, $d(Y,Z)=k-j$ and $d(X,Z)=k$ where
$1\leqslant j<k\leqslant n$. Then $\sigma(X,Y,Z)=1$.
\end{lemma}

\begin{pf}
Choose a hyperbolic basis $e_1,e_2,\ldots,e_n,f_1,f_2,\ldots,f_n$ for $V$, so that $B(e_i,e_j)=B(f_i,f_j)=0$ and $B(e_i,f_j)=\delta_{ij}$.
Since Sp$(2n,q)$ is transitive on triples of generators satisfying the given distance constraints, by Theorem~\ref{theorem4.4} we may suppose that
\[\begin{array}{c}
X=\langle e_1,e_2,\ldots,e_n\rangle,\quad Y=\langle f_1,f_2,\ldots,f_j,e_{j+1},e_{j+2},\ldots,e_n\rangle,\\
Z=\langle f_1,f_2,\ldots,f_k,e_{k+1},e_{k+2},\ldots,e_n\rangle.
\end{array}\]
We choose the determinant function $\delta_X$ on $X$ given by
\[\delta_X(x_1,x_2,\ldots,x_n)=\det\bigl[B(x_i,f_j)\,:\,1\leqslant i,j\leqslant n\bigr]\quad\hbox{for\ }x_1,x_2,\ldots,x_n\in X;\]
this is nothing other than the determinant of the $n\times n$ matrix whose columns are the coordinates of $x_1,\ldots,x_n$ with respect to the
basis $e_1,e_2,\ldots,e_n$. The determinant functions $\delta_Y$, $\delta_Z$ on $Y$ and on $Z$ are defined similarly, using the bases for $Y$ and on
$Z$ listed above. The computation of $\sigma(X,Z)$ is simplified by the fact that a basis for $X\cap Z$ is $e_{k+1},e_{k+2},\ldots,e_n$. We have
\[\delta_X(e_1,e_2,\ldots,e_n)=\delta_Z(f_1,\ldots,f_k,e_{k+1},\ldots,e_n)=1\]
and $\bigl[B(e_i,f_j)\,:\,1\leqslant i,j\leqslant k\bigr]$ is a $k\times k$ identity matrix, with determinant~1; thus $\sigma(X,Z)=1$. Exactly the same reasoning
gives $\sigma(X,Y)=\sigma(Y,Z)=1$, so $\sigma(X,Y,Z)=1$.\qed
\end{pf}

In the case of triples $X,Y,Z$ not lying on geodesic paths, however, $\sigma$ (or its two-graph $\Delta_\sigma$) yields interesting nontrivial information. In particular, the restriction of $\Delta_\sigma$
to partial spreads (sets of vertices of $\Gamma$ mutually at distance $n$) was investigated in~\cite[\S6]{M}.
Here we consider triangles in $\Gamma$:

\begin{lemma}\label{lemma4.6}
Suppose $q\equiv1\,$ \hbox{mod $4$}. Let $X,Y\in\V$ such that $d(X,Y)=1$, i.e.\ $X$ and $Y$ are adjacent in $\Gamma$.
There are $a_1=q{-}1$ common neighbors $Z$ of $X$ and $Y$ in $\Gamma$; and exactly half of the resulting triples $(X,Y,Z)$ are coherent.
\end{lemma}

\begin{pf}
Choose a hyperbolic basis $e_i,f_i$ as in the proof of Lemma~\ref{lemma4.5}. Again without loss of generality,
\[X=\langle e_1,e_2,\ldots,e_n\rangle,\quad Y=\langle f_1,e_2,\ldots,e_n\rangle,\quad Z=\langle e_1{+}\alpha f_1,e_2,\ldots,e_n\rangle\]
where $0\neq\alpha\in\Fq$. The $q{-}1$ choices of $\alpha$ give exactly the $q{-}1$ common neighbors of $X$ and $Y$ in $\Gamma$.
These bases of $X,Y,Z$ give rise to natural choices of determinant functions $\delta_X$,
$\delta_Y$, $\delta_Z$ as described in the proof of Lemma~\ref{lemma4.5}. When computing
$\sigma(X,Y),\sigma(Y,Z),\sigma(Z,X)$, we use $e_2,e_3,\ldots,e_n$ as the basis of $X\cap Y=X\cap Z=Y\cap Z$. Now
\[\delta_X(e_1,e_2,\ldots,e_n)=\delta_Z(e_1{+}\alpha f_1,e_2,\ldots,e_n)=1\]
and $B(e_1,e_1{+}\alpha f_1)=\alpha$, so $\sigma(X,Z)=\chi(\alpha)$. Similarly, $\sigma(X,Y)=\sigma(Y,Z)=1$ and
\[\sigma(X,Y,Z)=\chi(\alpha).\]
Since exactly half the nonzero elements of $\Fq$ are squares, the result follows.\qed
\end{pf}

\begin{theorem}\label{theorem4.7}
Suppose $q\equiv1\,$ \hbox{mod $4$}. Let $X,Y\in\V$ such that $d(X,Y)=k$. Then $Y$ has exactly $a_k=q^k{-}1$ neighbors $Z\in\V$
at distance $k$ from $X$ in $\Gamma$; and exactly half of the resulting triples $(X,Y,Z)$ are coherent.
\end{theorem}

\begin{pf}
The result holds for $k=1$ by Lemma~\ref{lemma4.6}, so we may assume $k\geqslant2$.
Given $X,Y\in\V$ with $d(X,Y)=k$, there are $\gauss{k}1$ choices of hyperplane $H<Y$ containing $X\cap Y$. Each such $H$ yields an Sp$(2,q)$-space
$H^\perp/H$, which contains $q{+}1$ subspaces of the form $Z/H$ with $Z\in\V$. One such $Z$ has distance $k{-}1$ from $X$, this being the subspace
$W=(Y{\cap}Z)+X\cap(Y{+}Z)=(Y{+}Z)\cap(X{+}(Y{\cap}Z))$. If we exclude $W$ and $Y$ itself, this leaves exactly $q{-}1$ choices of $Z$ having the required
distances from $X$ and $Y$; and this gives $(q{-}1)\gauss{k}1=q^k{-}1=a_k$ choices of $Z$, the full number. But for how many such $Z$ is the resulting
triple $(X,Y,Z)$ coherent? In each case $\sigma(X,W,Y)=\sigma(X,W,Z)=1$ by Lemma~\ref{lemma4.5}; therefore $\sigma(X,Y,Z)=\sigma(W,Y,Z)$.
But by Lemma~\ref{lemma4.6}, given
$W,Y$ at distance~1, exactly half of the $q{-}1$ choices of $Z$ yield coherent triples $(W,Y,Z)$. Therefore among the $a_k=(q{-}1)\gauss{k}1$ triples
$(X,Y,Z)$ with fixed $X$ and $Y$, exactly $\frac{q{-}1}{2}\gauss{k}1=(q^k{-}1)/2$ such triples are coherent.\qed
\end{pf}

\section{The Double Cover $\hatGamma\to\Gamma$}
\label{double_cover}

The resulting double cover $\hatGamma=\hatGamma(2n,q)\to\Gamma(2n,q)$ has vertex set $\hatV=\V\times\{\pm1\}$ and adjacency relation
\[(X,\varepsilon)\sim(Y,\varepsilon')\qquad\iff\qquad d(X,Y)=1\hbox{\ and\ }\varepsilon\varepsilon'=\sigma(X,Y).\]
The covering map is given by $(X,\varepsilon)\mapsto X$. 

\begin{theorem} \label{theorem5.1}
Every geodesic path
\[X_0\sim X_1\sim\cdots\sim X_k\]
in $\Gamma$ (meaning that $d(X_i,X_j)=|j-i|$) lifts to exactly two paths
\[(X_0,\varepsilon_0)\sim(X_1,\varepsilon_1)\sim\cdots\sim(X_k,\varepsilon_k)\]
in $\hatGamma$, in which $\varepsilon_k=\varepsilon_0\sigma(X_0,X_k)$ for each $k\geqslant1$; thus any one of the $\varepsilon_i$ determines all the others
along this path.
\end{theorem}

\begin{pf}
We have $\varepsilon_1=\varepsilon_0\sigma(X_0,X_1)$ by definition of adjacency in $\hatGamma$. Assuming that
$\varepsilon_i=\varepsilon_0\sigma(X_0,X_i)$ for some $i\in\{1,2,\ldots,k{-}1\}$,
\[\varepsilon_{i+1}=\varepsilon_i\sigma(X_i,X_{i+1})=\varepsilon_0\sigma(X_0,X_i)\sigma(X_i,X_{i+1})=\varepsilon_0\sigma(X_0,X_{i+1})\]
since $\sigma(X_0,X_i,X_{i+1})=1$ by Lemma~\ref{lemma4.5}.\qed
\end{pf}

\noindent However, not every geodesic path in $\hatGamma$ is obtained by lifting a geodesic path in~$\Gamma$. For example if $(X,Y,Z)$ is an incoherent
triangle in $\Gamma$, say with $\sigma(X,Y)=\varepsilon$, $\sigma(Y,Z)=\varepsilon'$ and $\sigma(X,Z)=-\varepsilon\varepsilon'$, then
\[(X,1)\sim(Y,\varepsilon)\sim(Z,\varepsilon\varepsilon')\sim(X,-1)\]
is a geodesic path of length~3 in $\hatGamma$, obtained by lifting a closed path of length~3 (not a geodesic path) in~$\Gamma$.

\begin{lemma} \label{lemma5.2}
Let $X_0\sim X_1\sim\cdots\sim X_k$ be a geodesic path of length $k\geqslant1$ in $\Gamma$, so that $d(X_i,X_j)=|j-i|$, and let
$\varepsilon,\varepsilon'\in\{\pm1\}$. Then $(X_0,\varepsilon)$ and $(X_k,\varepsilon')$ have distance $k$ or $k{+}1$ in $\hatGamma$, according
as $\sigma(X_0,X_k)=\varepsilon\varepsilon'$ or $-\varepsilon\varepsilon'$. In particular, the diameter of $\hatGamma$ is $\max\{n{+}1,3\}$.
\end{lemma}

\begin{pf}
If $\sigma(X_0,X_k)=\varepsilon\varepsilon'$, then we have a path
\[(X_0,\varepsilon_0)\sim(X_1,\varepsilon_1)\sim\cdots\sim(X_k,\varepsilon_k)\]
in $\hatGamma$ where $\varepsilon_i=\varepsilon_0\sigma(X_0,X_i)$ for $i=1,2,\ldots,k$; in particular if $\varepsilon_0=\varepsilon$ then
$\varepsilon_k=\varepsilon'$. Clearly this path in $\hatGamma$ is shortest possible.

Now suppose $\sigma(X_0,X_k)=-\varepsilon\varepsilon'$. We first obtain a path
\[(X_0,\varepsilon)\sim(X_1,\varepsilon_1)\sim\cdots\sim(X_{k-1},\varepsilon_{k-1})\]
in $\hatGamma$ where $\varepsilon_i=\varepsilon_0\sigma(X_0,X_i)$ for $i=1,2,\ldots,k{-}1$. Let $Y\in\V$ be adjacent to both $X_{k-1}$ and $X_k$ in
$\Gamma$, such that $\sigma(X_{k-1},Y,X_k)=-1$. (By Lemma~\ref{lemma4.6}, there are $\frac{q-1}{2}\geqslant1$ choices of such $Y\in\V$.) Appending
\[(X_{k-1},\varepsilon_{k-1})\sim(Y,\varepsilon'')\sim(X_k,\varepsilon'),\]
where $\varepsilon''=\varepsilon_{k-1}\sigma(X_{k-1},Y)=\varepsilon'\sigma(Y,X_k)$, we obtain a path of length $k{+}1$ from $(X_0,\varepsilon)$ to
$(X_k,\varepsilon')$ in $\hatGamma$; once again this path is shortest possible.\qed
\end{pf}

The fibers of the covering map $\hatGamma\to\Gamma$ are the {\it antipodal\/} pairs
$\{(X,1),(X,-1)\}$ for $X\in\V$.

\begin{lemma}\label{lemma5.3}
Let $(X,\varepsilon)$ and $(W,\varepsilon')$ be any two vertices of $\hatGamma$. Then
$(X,\varepsilon)$ and $(W,\varepsilon')$ are antipodal iff they are at distance~3 in $\hatGamma$ and are joined by exactly $\frac12 q(q^n{-}1)$ paths of length~3 in $\hatGamma$.
\end{lemma}

\begin{pf}Consider a typical antipodal pair $\{(X,1),(X,-1)\}$ where $X\in\V$. There exist $b_0=q\gauss{n}1$ vertices $Y\in\V$ adjacent to $X$ in $\Gamma$; and each such vertex $Y$ has $a_1=q{-}1$ neighbors $Z$ in common with $X$.
By Lemma~\ref{lemma4.6}, exactly half of these choices of the vertex $Z$ yield coherent triples $(X,Y,Z)$. In particular, $X$ lies in exactly $b_0{\cdot}\frac12 a_1=\frac12 q\bigl(q^n-1)$ incoherent triples $(X,Y,Z)$, giving
the same number of paths $(X,1)\sim(Y,\varepsilon)\sim(Z,\varepsilon')\sim(X,-1)$
of length~3 in $\hatGamma$. There is no path of length${}<3$ from $(X,1)$ to $(X,-1)$ in
$\hatGamma$, otherwise the covering map would give a closed
path of length${}<3$ from $X$ to $X$ in $\Gamma$. This shows that any two antipodal
vertices $(X,1),(X,-1)$ are at distance~3 in $\hatGamma$; and in each case there are exactly $\frac12 q(q^n{-}1)$ geodesic paths from $(X,1)$ to $(X,-1)$.

Conversely, let $(X,\varepsilon)$ and $(W,\varepsilon')$ be any two vertices at distance~3 in
$\hatGamma$. By Lemma~\ref{lemma5.2}, $d(X,W)\in\{0,2,3\}$
in $\Gamma$. Consider first the case that $d(X,W)=3$; then by Theorem~\ref{theorem5.1}, every geodesic path from $(X,\varepsilon)$ to
$(W,\varepsilon')$ in $\hatGamma$ arises from a unique geodesic path $X\sim Y\sim Z\sim W$ in $\Gamma$. There are exactly $c_3 c_2 c_1=(q^2{+}q{+}1)(q{+}1)$ such
geodesic paths from $X$ to $W$; and this number clearly cannot equal $\frac12 q(q^n-1)$.

Next suppose $d(X,W)=2$ in $\Gamma$. For every geodesic path
\[(X,\varepsilon)\sim(Y,\varepsilon'')\sim(Z,\varepsilon''')\sim(W,\varepsilon')\]
in $\hatGamma$, we have $X\sim Y\sim Z\sim W$ in $\Gamma$. Further, the condition $d(X,W)=2$ requires either $X\sim Z$ or $Y\sim W$ (but {\it not both,\/} by Lemma~\ref{lemma4.2}). We first count
geodesic paths satisfying $X\sim Z$, noting that the vertex $W$ has $c_2=q{+}1$ neighbors $Z$ in common with $X$; and in each case $\sigma(X,Z,W)=1$ by Lemma~\ref{lemma4.5}.
Moreover $Z$ has $a_1=q{-}1$ neighbors $Y$ in common with $X$ (all of which satisfy $\sigma(Y,Z,W)=1$, again by Lemma~\ref{lemma4.5}). By the two-graph
condition, we have $\sigma(X,Y,Z)=-1$ iff $\sigma(X,Y,W)=-1$. By Lemma~\ref{lemma4.6}, for each $Z$ there are exactly $\frac12\bigl(q{-}1\bigr)$
choices of $Y$ satisfying the latter condition; and each such pair $(Y,Z)$ yields a unique geodesic path
$(X,\varepsilon)\sim(Y,\varepsilon'')\sim(Z,\varepsilon''')\sim(W,\varepsilon')$. We obtain
$(q{+}1){\cdot}\frac12\bigl(q{-}1\bigr)=\frac12\bigl(q^2{-}1\bigr)$ geodesic paths in this case. There are another $\frac12\bigl(q^2{-}1\bigr)$ geodesic paths from $(X,\varepsilon)$ to $(W,\varepsilon')$ satisfying $Y\sim W$, for a total of $q^2{-}1$ geodesic paths. Once again, this number cannot equal $\frac12 q(q^n{-}1)$.\qed
\end{pf}

\begin{lemma}\label{lemma5.4}
$\Aut\hatGamma$ acts naturally on $\Gamma$, with kernel $\langle\zeta\rangle$,
inducing a {\em proper\/} subgroup $\Aut\hatGamma/\langle\zeta\rangle<\Aut\Gamma$.
\end{lemma}

\begin{pf}
By Lemma~\ref{lemma5.3}, $\Aut\hatGamma$ permutes fibres of the covering
map $\hatGamma\to\Gamma$, and so $\Aut\hatGamma$ acts naturally on $\Gamma$.
It remains to be shown that the induced subgroup $\Aut\hatGamma/\langle\zeta\rangle\leqslant\Aut\Gamma$ is proper.

Choose a hyperbolic basis $e_1,e_2,\ldots,e_n,f_1,f_2,\ldots,f_n$ for $V$, so that
$B(e_i,e_j)=B(f_i,f_j)=0$ and $B(e_i,f_j)=\delta_{ij}$, and let $\eta\in\Fq$ be a nonsquare.
Consider the subspaces $X,Y,Z,Z'\in\V$ defined by $X=\langle e_1,e_2,\ldots,e_n\rangle$,
$Y=\langle f_1,e_2,\ldots,e_n\rangle$, $Z=\langle e_1{+}f_1,e_2,\ldots,e_n\rangle$ and
$Z'=\langle e_1{+}\eta f_1,e_2,\ldots,e_n\rangle$. By straightforward computation,
$\sigma(X,Y,Z)=1$ and $\sigma(X,Y,Z')=-1$. Now consider $g\in GL(V)$ mapping our
original ordered basis to the new ordered basis $e_1,\eta e_2,\ldots,\eta e_n,\eta f_1,f_2,\ldots,f_n$ so that $B(u^g,v^g)=\eta B(u,v)$ for all $u,v\in V$; thus $g\in GSp(2n,q)$
induces an automorphism of the dual polar graph $\Gamma=\Gamma(2n,q)$.
However, $g$ maps the coherent triple $\{X,Y,Z\}$ to the non-coherent triple $\{X,Y,Z'\}$
and so does not preserve $\Delta_\sigma$. If $g$ were induced by an automorphism of
$\hatGamma$, this automorphism would map $\{X,Y,Z\}\times\{\pm1\}$ to
$\{X,Y,Z'\}\times\{\pm1\}$. However, the induced subgraphs of $\hatGamma$ on these
two 6-sets of vertices are not isomorphic (a pair of triangles and a 6-cycle, respectively;
see Section~\ref{graphs}).\qed
\end{pf}

The natural action of $\Sigma Sp(2n,q)$ on $\V$ lifts to an action on $\hatV$ as follows:
Let $g\in\Sigma Sp(2n,q)$ with associated field automorphism $\tau_g$ in the earlier notation of this section.
Given $(U,\varepsilon)\in\hatV$, the map
\[U^n\to\Fq,\quad (u_1,u_2,\ldots,u_n)\mapsto\delta_{U^g}(u_1^g,u_2^g,\ldots,u_n^g)^{\tau_g^{-1}}\]
is a determinant function; so there exists a nonzero constant $\lambda_{g,U}\in\Fq$ such that
\[\delta_{U^g}(u_1^g,u_2^g,\ldots,u_n^g)=\lambda_{g,U}\delta_U(u_1,u_2,\ldots,u_n)^{\tau_g}.\]
Define $(U,\varepsilon)^g=(U^g,\chi(\lambda_{g,U})\varepsilon)$. One easily checks that
this defines an action of $\Sigma Sp(2n,q)$
on $\hatV$. The central element $-I\in Sp(2n,q)$ fixes every $U\in\V$ and since
\[\delta_U(-u_1,-u_2,\ldots,-u_n)=(-1)^n\delta_U(u_1,u_2,\ldots,u_n)\]
where $\chi(-1)^n=1$,
$-I$ acts trivially on $\hatV$; thus $\Sigma Sp(2n,q)$ induces a permutation group
$P\Sigma Sp(2n,q)$ on~$\hatV$. The transposition $\zeta$ which exchanges antipodal
vertices via $(U,1)\stackrel\zeta\leftrightarrow(U,-1)$ is not induced by any element of
$P\Sigma Sp(2n,q)$ since $Z(P\Sigma Sp(2n,q))=1$, so we obtain a permutation group
$\langle\zeta\rangle\times P\Sigma Sp(2n,q)$ acting faithfully on $\hatV$. We show that this
permutation group preserves the graph $\hatGamma$, and is in fact its full automorphism group:

\begin{theorem}\label{theorem5.5}
$\Aut\hatGamma\cong2\times P\Sigma Sp(2n,q)$ where this group acts as defined above.
The full automorphism group of the two-graph associated to $\sigma$ is $\Aut\Delta_\sigma\cong P\Sigma Sp(2n,q)$.
\end{theorem}

\begin{pf}
Suppose $(X,\varepsilon)\sim(Y,\varepsilon')$ in $\hatGamma$, so that
$\sigma(X,Y)=\varepsilon\varepsilon'$; and let $g\in\Sigma Sp(2n,q)$ with
$\tau_g\in\Aut\Fq$ as above. Then by Proposition~\ref{prop4.3}(iv) we have
\[\sigma(X^g,Y^g)=\chi(\lambda_{g,X}\lambda_{g,Y})\sigma(X,Y)
=\bigl(\chi(\lambda_{g,X})\varepsilon\bigr)\bigl(\chi(\lambda_{g,Y})\varepsilon'\bigr)\]
so that by definition, $(X,\varepsilon)^g\sim(Y,\varepsilon')^g$. Thus $P\Sigma Sp(2n,q)$,
acting on $\hatGamma$ as defined above, preserves the graph $\hatGamma$.
It is clear that the central factor $\langle\zeta\rangle$ also preserves $\hatGamma$,
so that $\Aut\hatGamma$ has a subgroup isomorphic to
$\langle\zeta\rangle\times P\Sigma Sp(2n,q)$. Moreover by Proposition~\ref{prop3.1},
$\Aut\hatGamma/\langle\zeta\rangle\cong\Aut\Delta_\sigma$. (We use the fact that by
Lemma~\ref{lemma5.4}, $\Aut_\zeta\hatGamma=\Aut\hatGamma$ in the notation
of Proposition~\ref{prop3.1}.)

Suppose now that $n\geqslant2$, so that $\Aut\Gamma\cong P\Gamma Sp(2n,q)$ by
Theorem~\ref{theorem4.1}. By Lemma~\ref{lemma5.4}, $\Aut\hatGamma$
acts on $\Gamma$, inducing a group of automorphisms satisfying
\[P\Sigma Sp(2n,q)\leqslant\Aut\hatGamma/\langle\zeta\rangle<P\Gamma Sp(2n,q).\]
This forces $\Aut\hatGamma\cong\langle\zeta\rangle\times P\Sigma Sp(2n,q)$ and
$\Aut\Delta_\sigma\cong P\Sigma Sp(2n,q)$.

Finally suppose $n=1$, so that $\Delta_\sigma$ is the
Taylor-Paley two-graph on $q+1$ vertices, with full automorphism group
$\Aut\Delta_\sigma\cong P\Sigma Sp(2,q)=P\Sigma L(2,q)$ by \cite[Theorem 2]{Ta};
see also~\cite[\S4]{M}. As above, $\Aut\hatGamma$ has a subgroup
isomorphic to $\langle\zeta\rangle\times P\Sigma Sp(2,q)$, and
$\Aut\hatGamma/\langle\zeta\rangle\cong\Aut\Delta_\sigma\cong P\Sigma Sp(2,q)$,
so we must have equality:
$\Aut\hatGamma\cong2\times P\Sigma Sp(2,q)=2\times P\Sigma L(2,q)$.\qed
\end{pf}

\section{The Association Scheme}
\label{scheme}

From the double cover $\hatGamma\to\Gamma$, we now construct association schemes.
As we will see in Section~\ref{Q-poly}, this gives
a new family of $Q$-polynomial association schemes. We begin with the
relevant definitions, following \cite[Chapter~2]{BCN}.

Let $\Omega$ be a finite set. A {\em (symmetric) $d$-class association scheme on $\Omega$\/} is a pair $(\Omega,\Rel)$ such that
\begin{enumerate}
\item $\Rel=\{R_0,\ldots,R_d\}$ is a partition of $\Omega\times\Omega$;
\item $R_0$ is the identity relation on $\Omega$;
\item $R_i=R_i^{\top}$ for $0\leqslant i\leqslant d$; and
\item there are constants $p_{ij}^k$ such that for any pair $(x,y)\in R_k$, the number of $z \in\Omega$ such that $(x,z)\in R_i$ and $(z,y)\in R_j$ equals $p_{ij}^k$.
\end{enumerate}
For the rest of this paper, all association schemes are symmetric (i.e.\ the third property above holds).
Each relation $R_i$ has adjacency matrix $A_i$ defined by
\[(A_i)_{x,y}=\begin{cases}
1,&\mbox{if\ }(x,y)\in R_i;\\
0,&\mbox{otherwise.}
\end{cases}\]
The axioms above imply that $p_{ji}^k=p_{ij}^k$ and the matrices $A_0,\ldots,A_d$ form an algebra of symmetric matrices satisfying
$A_iA_j=\sum_k p_{ij}^k A_k$. This matrix algebra is also closed under Schur (entrywise) multiplication, which we will denote by 
`$\circ$'. This algebra is referred to as the {\em Bose-Mesner algebra\/} $\BM$ of the association scheme.  

Since $\BM$ is a commutative algebra consisting of symmetric matrices, its elements are simultaneously diagonalizable, and $\BM$ has a second basis consisting of primitive idempotents $E_0,\ldots,E_d$.  
We define the parameters $Q_{ij}$ by $E_j =\frac{1}{|\Omega|}\sum_i Q_{ij}A_i$.  Similarly we define the parameters $P_{ij}$ by the relation 
$A_j=\sum_i P_{ij} E_i$.  The matrix $P$ of parameters $P_{ij}$ is often referred to as the character table of the scheme.   
The matrix $Q$ of parameters $Q_{ij}$ satisfies $Q=|\Omega|P^{-1}$. 

We say an association scheme is {\em$Q$-polynomial\/} if, after suitably reindexing its idempotents, 
the idempotent $E_j$ is a degree $j$ polynomial in $E_1$ (where multiplication is done entrywise). This is equivalent to the condition
that the $j$th column of the $Q$-matrix is a degree $j$ polynomial of the column 1 of the $Q$-matrix (note that we start indexing the columns at 0).

Permutation groups give many examples of association schemes.  Let $G$ be a transitive
permutation group acting on a finite set $\Omega$, and suppose the orbits of $G$ on
$\Omega\times\Omega$ happen to be symmetric relations; such a group is called {\em generously transitive\/}).
It is not hard to check that the orbits of $G$ on $\Omega\times\Omega$ form an association scheme. We will refer to these schemes as {\em Schurian schemes\/}.
\bigskip

We will now construct a $(2n{+}1)$-class Schurian association scheme $\mathcal{S}=\mathcal{S}_{n,q}$ with vertex set $\hatV=\V\times\{\pm1\}$ of cardinality
$|\hatV|=2|\V|=2q^{n^2}\!\prod_{i=1}^n(q^{2i}-1)$ using the Maslov index $\sigma$
and the double cover $\hatGamma\to\Gamma$ (defined as in Section~\ref{double_cover}).
For $k=0,1,2,\ldots,n$, the $k$th and $(2n{+}1{-}k)$th relations are given by
\begin{align*}
R_k&=\{((X,\varepsilon),(Y,\varepsilon'))\in\hatV\times\hatV\;:\;d(X,Y)=k,\;\varepsilon\varepsilon'=\sigma(X,Y)\};\\
R_{2n+1-k}&=\{((X,\varepsilon),(Y,\varepsilon'))\in\hatV\times\hatV\;:\;d(X,Y)=k,\;\varepsilon\varepsilon'=-\sigma(X,Y)\}.
\end{align*}
These are symmetric relations which clearly partition $\hatV\times\hatV$. In particular,
$R_1$ is the adjacency relation of our graph $\hatGamma$ of Section~\ref{double_cover}; and the identity and antipodality relations are
\begin{align*}
R_0&=\{((X,\varepsilon),(X,\varepsilon)):X\in\V,\;\varepsilon=\pm1\};\\
R_{2n+1}&=\{((X,\varepsilon),(X,-\varepsilon)):X\in\V,\;\varepsilon=\pm1\}.
\end{align*}
We will write
\[(X,\varepsilon)\irel{i}(Y,\varepsilon')\iff((X,\varepsilon),(Y,\varepsilon'))\in R_i\,.\]
In the following, the parameters $a_i,b_i,c_i$ are those of the dual polar graph $\Gamma$
as given in Section~\ref{dual_polar}.

\begin{lemma}\label{lemma6.1}
Let $(X,\varepsilon)\irel{k}(Y,\varepsilon')$ where $k\in\{0,1,2,\ldots,2n{+}1\}$. The number of $(Z,\varepsilon'')\in\hatV$
such that $(X,\varepsilon)\irel{i}(Z,\varepsilon'')\irel{1}(Y,\varepsilon')$ is
\[p_{i,1}^k=\begin{cases}
c_k=\gauss{k}1,&\mbox{if }\,i=k{-}1\leqslant n;\\
\noalign{\smallskip}
\frac12 a_k=\frac{1}{2}\bigl(q^k-1),&\mbox{if }\,i=k;\\
\noalign{\smallskip}
b_k=q^{k+1}\gauss{n-k}1,&\mbox{if }\,i=k{+}1\leqslant n{+}2;\\
\noalign{\smallskip}
b_{2n+1-k}=q^{2n+2-k}\gauss{k-n-1}1,&\mbox{if }i=k{-}1\geqslant n{-}1;\\
\noalign{\smallskip}
\frac12 a_{2n+1-k}=\frac{1}{2}\bigl(q^{2n+1-k}-1),&\mbox{if }\,i=2n{+}1{-}k;\\
\noalign{\smallskip}
c_{2n+1-k}=\gauss{2n+1-k}1,&\mbox{if }\,i=k{+}1\geqslant n{+}1;\\
\noalign{\smallskip}
0,&\hbox{otherwise.}
\end{cases}\]
\end{lemma}

\begin{pf}(i) First suppose $d(X,Y)=k\leqslant n$, so $\varepsilon\varepsilon'=\sigma(X,Y)$. Then $(Z,\varepsilon'')\in\hatV$ satisfies
$(X,\varepsilon)\irel{i}(Z,\varepsilon'')\irel{1}(Y,\varepsilon')$ iff
\[\begin{cases}
\mbox{\qquad case (i.a)}\\
\hline
i=d(X,Z)\leqslant n\\
d(Z,Y)=1\\
\varepsilon''=\varepsilon\sigma(X,Z)=\varepsilon'\sigma(Y,Z)
\end{cases}\mbox{\quad\em or\quad}
\begin{cases}
\mbox{\qquad case (i.b)}\\
\hline
i=2n{+}1{-}d(X,Z)\geqslant n{+}1\\
d(Z,Y)=1\\
\varepsilon''=-\varepsilon\sigma(X,Z)=\varepsilon'\sigma(Y,Z)
\end{cases}\]
Moreover, each such $(Z,\varepsilon'')$ satisfies $d(X,Z)\in\{k{-}1,k,k{+}1\}$ by the triangle inequality.

There are exactly $c_k=\gauss{k}1$ choices of $Z\in\V$ satisfying $d(X,Z)=k{-}1$ and $d(Z,Y)=1$. Each such $Z$ yields a coherent
triple $(X,Y,Z)$ by Lemma~\ref{lemma4.5}, so $\varepsilon\varepsilon'=\sigma(X,Y)=\sigma(X,Z)\sigma(Y,Z)$. This yields $c_k$ pairs $(Z,\varepsilon'')$,
all of which satisfy~(i.a).

There are exactly $b_k=q^{k+1}\gauss{n-k}1$ choices of $Z\in\V$ satisfying $d(X,Z)=k{+}1$ and $d(Z,Y)=1$. Each such $Z$ yields a coherent
triple $(X,Y,Z)$ by Lemma~\ref{lemma4.5}, once again with $\varepsilon\varepsilon'=\sigma(X,Y)=\sigma(X,Z)\sigma(Y,Z)$. This yields $b_k$ pairs $(Z,\varepsilon'')$,
all of which satisfy~(i.a).

There are exactly $a_k=q^k-1$ choices of $Z\in\V$ satisfying $d(X,Z)=k$ and $d(Z,Y)=1$. By Theorem~\ref{theorem4.7}, exactly $a_k/2$ of these $Z$ yield coherent
triples $(X,Y,Z)$, in which case $\varepsilon\varepsilon'=\sigma(X,Y)=\sigma(X,Z)\sigma(Y,Z)$; this yields $a_k/2$ pairs $(Z,\varepsilon'')$
satisfying~(i.a). The remaining $a_k/2$ of these $Z$ yield incoherent
triples $(X,Y,Z)$, with $\varepsilon\varepsilon'=\sigma(X,Y)=-\sigma(X,Z)\sigma(Y,Z)$; and the resulting pairs $(Z,\varepsilon'')$ satisfy~(i.b).\qed
\end{pf}

In Section~\ref{double_cover} we lifted the action of $P\Sigma Sp(2n,q)$ on $\V$, to a
transitive permutation action of $\langle\zeta\rangle\times P\Sigma Sp(2n,q)$ on $\hatV$
(below Lemma~\ref{lemma5.4}). Theorem~\ref{theorem5.5} shows
that this group preserves $R_1$ (the adjacency relation of the graph $\hatGamma$).
We next show that this group preserves {\em each\/} of the relations $R_i$, and
so gives the full automorphism group of the scheme.

\begin{lemma}\label{lemma6.2}
The diagonal action of $2\times PSp(2n,q)$ on $\hatV\times\hatV$ preserves each
of the relations $R_i$. The same conclusion holds for the subgroup $2\times P\Sigma Sp(2n,q)$.
\end{lemma}

\begin{pf}
Clearly the central factor $(U,\varepsilon)\stackrel\zeta\leftrightarrow(U,-\varepsilon)$ preserves each $R_i$.
Now let $g\in Sp(2n,q)$, and suppose
$X,Y\in\V$ such that $d(X,Y)=k\in\{0,1,2,\ldots,n\}$. Also let $\varepsilon,\varepsilon'\in\{\pm1\}$, so that $((X,\varepsilon),(Y,\varepsilon'))\in R_k$ or $R_{2n+1-k}$
according as $\varepsilon\varepsilon'\sigma(X,Y)=1$ or $-1$.
Since $g$ preserves distances in $\Gamma$, $d(X^g,Y^g)=k$. Let $x_1,x_2,\ldots,x_n$ and
$y_1,y_2,\ldots,y_n$ be bases for $X$ and $Y$ respectively, such that a basis for
$X\cap Y$ is formed by $x_{k+1}{=}y_{k+1}$, $x_{k+2}{=}y_{k+2}$, \dots,
$x_n{=}y_n$. Then $(X,\varepsilon)^g=(X^g,\chi(\lambda_{g,X})\varepsilon)$ and
$(Y,\varepsilon')^g=(Y^g,\chi(\lambda_{g,Y})\varepsilon')$ where
\begin{align*}
&\chi(\lambda_{g,X})\varepsilon\chi(\lambda_{g,Y})\varepsilon'\sigma(X^g,Y^g)\\
&\qquad=\varepsilon\varepsilon'\chi\bigl(\lambda_{g,X}\delta_{X^g}(x_1^g,\ldots,x_n^g)\lambda_{g,Y}\delta_{Y^g}(y_1^g,\ldots,y_n^g)\\
&\qquad\phantom{{}={}}\times\det\bigl[B(x_i^g,y_j^g):1\leqslant i,j\leqslant k\bigr]\bigr)\\
&\qquad=\varepsilon\varepsilon'\chi\bigl(\lambda_{g,X}^2\lambda_{g,Y}^2\delta_X(x_1,\ldots,x_n)\delta_Y(y_1,\ldots,y_n)\det\bigl[B(x_i,y_j)\,:\,1\leqslant i,j\leqslant k\bigr]\bigr)\\
&\qquad=\varepsilon\varepsilon'\sigma(X,Y).
\end{align*}
If this value is $1$, then both $(X,\varepsilon)\irel{k}(Y,\varepsilon')$ and
$(X,\varepsilon)^g\irel{k}(Y,\varepsilon')^g$; but if the latter value is $-1$, then
$(X,\varepsilon)\irel{\scriptscriptstyle{2n{+}1{-}k}}(Y,\varepsilon')$ and
$(X,\varepsilon)^g\irel{\scriptscriptstyle{2n{+}1{-}k}}(Y,\varepsilon')^g$.

Thus $2\times PSp(2n,q)$ preserves the relations $R_i$ as claimed. A similar argument
holds for $2\times P\Sigma Sp(2n,q)$.\qed
\end{pf}

It is easy to see that $\langle\zeta\rangle\times P\Sigma Sp(2n,q)$ acts transitively on
each $R_i$, and similarly for $\langle\zeta\rangle\times PSp(2n,q)$. This yields

\begin{theorem}\label{theorem6.3}
The diagonal action of the group $2\times PSp(2n,q)$ on $\hatV\times\hatV$ has orbits $R_0$, $R_1$, \dots, $R_{2n+1}$; so these form the relations of a $(2n+1)$-class Schurian association scheme. The same conclusion holds for $2\times P\Sigma Sp(2n,q)$,
which is therefore the full automorphism group of the association scheme $\mathcal{S}$.\qed
\end{theorem}

\section{The $Q$-polynomial property}
\label{Q-poly}

In this section we will use some parameters of the scheme to prove that the association scheme $\mathcal{S}$ is $Q$-polynomial.
We will benefit from the action of the $A_i$'s by left-multiplication on the Bose-Mesner algebra,
resulting in matrices $L_i$ defined by $(L_i)_{kj}=p_{ij}^k$.
In particular, the parameter $p_{1j}^k$ of the scheme from Lemma~\ref{lemma6.1}, is the $(k,j)$-entry of the matrix
\[L_1=\begin{pmatrix}  
0&b_0&0&0&0&\cdots &0&0&0&0 \\
1&\frac{a_1}{2}&b_1&0&0&\cdots &0&0&\frac{a_1}{2}&0 \\
0&c_2&\frac{a_2}{2}&\ddots&0&\cdots &0&\iddots&0&0 \\
0&0&\ddots&\ddots&b_{d-1}&0&\frac{a_{d-1}}{2}&0&0&0 \\
0&0&0&c_d&\frac{a_d}{2}&\frac{a_d}{2} &0&0&0&0 \\
0&0&0&0&\frac{a_d}{2}&\frac{a_d}{2} &c_d&0&0&0 \\
0&0&0&\frac{a_{d-1}}{2}&0&b_{d-1}&\ddots&\ddots&0&0 \\
0&0&\iddots&0&0&0 &\ddots&\frac{a_2}{2}&c_2&0 \\
0&\frac{a_1}{2}&0&0&0&0 &0&b_1&\frac{a_1}{2}&1 \\
0&0&0&0&0&0 &0&0&b_0&0
\end{pmatrix}.\]

As it turns out, this matrix has distinct eigenvalues, which in turn will give us a great deal of information about the scheme.
In particular by \cite[Proposition 2.2.2]{BCN}, the columns of $Q$ are right eigenvectors of $L_1$.
We will use the following generalization of \cite[Theorem 8.1.1]{BCN} to prove that $\mathcal{S}$ is $Q$-polynomial.

\begin{theorem}\label{theorem7.1}
Suppose  $A_i$ is a matrix in a $d$-class association scheme $(\Omega,\Rel)$ with $d+1$ distinct eigenvalues.
Then $(\Omega,\Rel)$ is $Q$-polynomial if and only if there is a sequence of {\em distinct\/} complex numbers $\sigma_0,\sigma_1,\ldots,\sigma_d$ and polynomials 
$s_0(x),s_1(x),\ldots,s_d(x)$ of degree $0,1,\ldots,d$, respectively, with
\[\sum_j p_{ij}^k\sigma_j^\ell=s_\ell(\sigma_k)\]
for $0 \leqslant\ell\leqslant d$. Furthermore, the leading coefficients of the polynomials $s_0(x)$, $s_1(x)$, \dots, $s_d(x)$ are precisely the eigenvalues of $A_i$ in a $Q$-polynomial ordering.
\end{theorem}

\begin{pf}
Without loss of generality we assume $A_1$ has this property.  Let $L_1$ be the corresponding intersection matrix.  Let 
$S,T$ be the $d+1$ by $d+1$ matrices with $S_{jk}=\sigma_j^k$ and $T_{jk}$ equal to the coefficient of $x^j$ in the polynomial $s_k(x)$. Then the above statement is equivalent to 
$L_1 S=ST$.  Then $L_i$ is similar to $T$ and since $T$ is upper triangular, the diagonal entries of $T$ are precisely the eigenvalues of $L_1$. 
Since $T$ is an upper triangular matrix with distinct diagonal entries, an easy induction shows that it can be diagonalized by an upper triangular matrix.  Namely, there is an invertible matrix $U$ and a diagonal matrix $D$ with $D_{jj}=T_{jj}$ such that $U^{-1}TU=D$.
Then $L_1(SU)=(SU)D$. This implies that the columns of $SU$ are eigenvectors of $L_1$, hence there is a diagonal matrix $D'$ such that $SUD'=Q$.
Since the $j$th column of $SUD'$ is a degree $j$ polynomial of the first column of $T$, which is a linear combination of columns 0 and 1 of $SUD'$, it is clear that the $j$th column of $SUD'$ is a degree $j$ polynomial of the 
first column of $SUD'$.  This implies that the columns of $Q$ are in a given $Q$-polynomial ordering, which in turn implies that the ordering of the eigenvalues in $T$ is a $Q$-polynomial ordering.\qed
\end{pf}

This leads to our main result:

\begin{theorem}\label{theorem7.2}
The scheme $\mathcal{S}$ is $Q$-polynomial. Furthermore, it has two $Q$-polynomial orderings.
\end{theorem}

\begin{pf}
Let $r=\sqrt{q}$ and $d=2n+1$. We define the sequence of polynomials
\[s_\ell(x)=\begin{cases}
r^\ell\gauss{n-\ell+1}{1}x^\ell+\frac{1}{r^{\ell-2}}\gauss{\ell-1}{1}x^{\ell-2},&\mbox{for $\ell$ odd};\\
\noalign{\smallskip}
r^\ell\left(\gauss{n-\ell+1}{1}-\frac{1}{r^\ell}\right)x^\ell+\frac{1}{r^{\ell-2}}\left( \gauss{\ell-1}{1}+r^{\ell-2} \right)x^{\ell-2},&\mbox{for $\ell$ even}
\end{cases}\]
and constants
\[\sigma_j=\begin{cases}
\frac{1}{r^j},&\mbox{for\ }0\leqslant j\leqslant n;\\
-\frac{1}{r^{2n+1-j}},&\mbox{for\ }n+1\leqslant j\leqslant2n+1.
\end{cases}\]
The polynomials $s_0(x),\ldots,s_{2n+1}(x)$ realize $\sigma_0,\ldots,\sigma_{2n+1}$
as a $Q$-sequence for $\mathcal{S}$, as we proceed to show by direct computation.
For $k\leqslant n$ we have $\sum_j p_{1j}^k \sigma_j^\ell=c_k \sigma_{k-1}^\ell+\frac{a_k}{2} \sigma_{k}^\ell+b_k \sigma_{k+1}^\ell+\frac{a_k}{2} \sigma_{2n+1-k}^\ell $.
For odd $\ell$ this reduces to
\begingroup
\allowdisplaybreaks
\begin{align*}
c_k \sigma_{k-1}^\ell+b_k \sigma_{k+1}^\ell
&=\tfrac{1}{r^{(k-1)\ell}}\gauss{k}{1}+q\left(\gauss{k}{1}-\gauss{n}{1}-\gauss{k}{1} \right)\tfrac{1}{r^{(k+1)\ell}}\\
&=\tfrac{1}{r^{(k+1)\ell}}\left(\gauss{k}{1}q^\ell+q\bigl(\gauss{n}{1}-\gauss{k}{1}\bigr)\right)\\
&=\tfrac{1}{r^{(k+1)\ell}}\left(\gauss{n+1}{1}-\gauss{\ell}{1}+\gauss{k+\ell}{1}- \gauss{k+1}{1}\right)\\
&=\tfrac{1}{r^{(k+1)\ell}}\left( \gauss{n-\ell+1}{1}r^{2\ell}+\gauss{\ell-1}{1}r^{2(k+1)}\right)\\
&=r^\ell\gauss{n-\ell+1}{1}\tfrac{1}{r^{k\ell}}+\tfrac{1}{r^{\ell-2}}\gauss{\ell-1}{1}\tfrac{1}{r^{k(\ell-2)}}\\
&=s_\ell(\sigma_k),
\end{align*}
whereas for even $\ell$ we have
\begin{align*}
\textstyle{\sum_j p_{1j}^k \sigma_j^\ell}
&=c_k\sigma_{k-1}^\ell+a_k \sigma_{k}^\ell+b_k \sigma_{k+1}^\ell\\
&=\gauss{k}{1}\tfrac{1}{r^{(k-1)\ell}}+(q-1)\gauss{k}{1}\tfrac{1}{r^{k\ell}}+q\left(\gauss{k}{1}-\gauss{n}{1}-\gauss{k}{1}\right)\tfrac{1}{r^{(k+1)\ell}}\\
&=\tfrac{1}{r^{(k+1)\ell}}\left(\gauss{k}{1}q^{\ell}+q\bigl(\gauss{n}{1}-\gauss{k}{1}\bigr)+r^{2k+\ell}-r^\ell\right)\\
&=\tfrac{1}{r^{(k+1)\ell}}\left(\gauss{n+1}{1}-\gauss{\ell}{1}+\gauss{k+\ell}{1}- \gauss{k+1}{1}+r^{2k+\ell}-r^\ell\right)\\
&=\tfrac{1}{r^{(k+1)\ell}}\left(\gauss{n-\ell+1}{1}r^{2\ell}-r^\ell+\gauss{\ell-1}{1}r^{2(k+1)}+r^{2k+\ell}\right)\\
&=r^\ell\left(\gauss{n-\ell+1}{1}-\tfrac{1}{r^\ell}\right)\tfrac{1}{r^{k\ell}}+\tfrac{1}{r^{\ell-2}}\left(\gauss{\ell-1}{1}+r^{\ell-2}\right)\tfrac{1}{r^{k(\ell-2)}}\\
&=s_\ell(\sigma_k).
\end{align*}
\endgroup

Now we deal with $k\geqslant n+1$, noting that
\[\sum_j p_{1j}^k \sigma_j^\ell=b_{2n+1-k} \sigma_{k-1}^\ell+\tfrac{a_{2n+1-k}}{2} \sigma_{k}^\ell+c_{2n+1-k} \sigma_{k+1}^\ell+\tfrac{a_{2n+1-k}}{2} \sigma_{2n+1-k}^\ell.\]
For odd $\ell$ this reduces to
\begin{align*}
b_{2n+1-k} \sigma_{k-1}^\ell+c_{2n+1-k} \sigma_{k+1}^\ell&=- b_{2n+1-k} \sigma_{2n+2-k}^\ell - c_{2n+1-k} \sigma_{2n-k}^\ell\\
&=-s_\ell( \sigma_{2n+1-k})=s_\ell(-\sigma_{2n+1-k})=s_\ell(\sigma_{k}),
\end{align*}
while for even $\ell$ we obtain
\begin{align*}
b_{2n+1-k}\sigma_{k-1}^\ell+{}&a_{2n+1-k}\sigma_{k}^\ell+c_{2n+1-k} \sigma_{k+1}^\ell\\
&=b_{2n+1-k}\sigma_{2n+2-k}^\ell+a_{2n+1-k}\sigma_{2n+1-k}^\ell+c_{2n+1-k}\sigma_{2n-k}^\ell\\
&=s_\ell( \sigma_{2n+1-k})=s_\ell(-\sigma_{2n+1-k})=s_\ell(\sigma_{k}).
\end{align*}

For nonsquare $q$ the splitting field of $\mathcal{S}$ is irrational, implying that it is a quadratic extension of the rationals, namely 
$\mathbb{Q}(r)$. The Galois group acts faithfully on the idempotents of the scheme, yielding a second $Q$-polynomial ordering.  
This second $Q$-polynomial ordering can also be obtained by replacing $r\mapsto-r$ in both the $\sigma_j$ and the polynomials $s_\ell(x)$, showing that this second ordering exists for square $q$ as well. \qed
\end{pf}

We note that by a result of Suzuki \cite{SUZ}, $Q$-polynomial schemes can have at most two $Q$-polynomial orderings.

\section{The $P$-matrix}
\label{P-matrix}

We now compute the $P$-matrix of the scheme $\mathcal{S}$, expressing it in terms of the auxiliary matrices $\Pone$ and $\Ptwo$ whose entries are defined by
\begin{align*}
\Pone_{ij}&=\sum\limits_{l=0}^j (-1)^\ell r^{j-2\ell+(j-\ell)^2+\ell^2}\gauss{i}{\ell}\gauss{n-i}{j-\ell};\\
\Ptwo_{ij}&=\sum\limits_{\ell=0}^j (-1)^\ell r^{(j-\ell)^2+\ell^2}\gauss{i}{\ell}\gauss{n-i}{j-\ell}.
\end{align*}
By \cite[Proposition 2.2.2]{BCN}, the $P$-matrix is determined by the left-normalized left eigenvectors of $L_1$. We first show that the rows of $\Pone$ and $\Ptwo$ are left
eigenvectors of the matrices defined by
\[\Mone=\begin{pmatrix}  
0&b_0&&&&\\
1&a_1&b_1&&&\\
&c_2&a_2&\smash\ddots&&\\
\noalign{\vskip-2pt}
&&\ddots&\ddots&\ddots&\\
\noalign{\vskip-3pt}
&&&\ddots&a_{d-1}&b_{d-1}\\
&&&&c_{d}&a_d
\end{pmatrix},\quad
\Mtwo=\begin{pmatrix}  
0&b_0&&&&\\
1&0&b_1&&&\\
\noalign{\vskip-4pt}
&c_2&0&\ddots&&\\
\noalign{\vskip-2pt}
&&\ddots&\ddots&\ddots&\\
\noalign{\vskip-5pt}
&&&\ddots&0&b_{d-1}\\
&&&&c_{d}&0
\end{pmatrix}\]
respectively. We will show that the corresponding diagonal forms are
\[\Done=\diag(\Pone_{i0},\Pone_{i1},\ldots,\Pone_{in}),\quad\Dtwo=\diag(\Ptwo_{i0},\Ptwo_{i1},\ldots,\Ptwo_{in}).\]
The ordering we give to the eigenvectors of $\Mone$ and $\Mtwo$ may seem arbitrary, but will be important later.

\begin{theorem}\label{theorem8.1}
$\Pone\Mone=\Done\Pone$ and $\Ptwo\Mtwo=\Dtwo\Ptwo$.
\end{theorem}

\begin{pf}
Fix $i$ and let $v_i=(\Pone_{i0},\Pone_{i1},\ldots,\Pone_{in})$.
We must show that $v_i\Mone=\Pone_{i1}v_i$. In particular, we need to show the following recurrence holds for all $j$:
\begin{align*}
b_{j-1}\sum\limits_{\ell=0}^{j-1}(-1)^\ell&r^{j-1-2\ell +(j-1-\ell)^2+\ell^2}
\gauss{i}{\ell}\gauss{n-i}{j-1-\ell} + a_j \sum\limits_{\ell=0}^{j}(-1)^\ell r^{j-2\ell +(j-\ell)^2+\ell^2}\gauss{i}{\ell}\gauss{n-i}{j-\ell}\\
&\qquad+c_{j+1}\sum\limits_{\ell=0}^{j+1}(-1)^\ell r^{j+1-2\ell+(j+1-\ell)^2+\ell^2}\gauss{i}{\ell}\gauss{n-i}{j+1-\ell}\\
&=\Pone_{i1} \sum\limits_{\ell=0}^j(-1)^\ell r^{j-2\ell+(j-\ell)^2+\ell^2}\gauss{i}{\ell}\gauss{n-i}{j-\ell}.
\end{align*}
Multiplying both sides by $q-1$ and substituting for $b_{j-1}, a_j$ and $c_{j+1}$, we find this is equivalent to showing that the quantity $z_j$, defined as follows, vanishes for all $j$:
\begin{align*}
z_j={}&(q^{n+1}-q^j)\sum\limits_{\ell=0}^{j-1}(-1)^\ell r^{j-2\ell-1+(j-\ell-1)^2+\ell^2} \gauss{i}{\ell}\gauss{n-i}{j-\ell-1}\\
&+(q-1)(q^j-1)\sum\limits_{\ell=0}^{j}(-1)^\ell r^{j-2\ell+(j-\ell)^2+\ell^2}\gauss{i}{\ell}\gauss{n-i}{j-\ell}\\
&+(q^{j+1}-1) \sum\limits_{\ell=0}^{j+1}(-1)^\ell r^{ j-2\ell+1+(j-\ell+1)^2+\ell^2}\gauss{i}{\ell}\gauss{n-i}{j-\ell+1}\\
&-\bigl(q^i(q^{n-2i+1}-1)-q+1\bigr)\sum\limits_{\ell=0}^j (-1)^\ell r^{ j-2\ell+(j-\ell)^2+\ell^2} \gauss{i}{\ell}\gauss{n-i}{j-\ell}.
\end{align*}
The second and last sums combine, simplifying to
\begin{align*}
z_j={}&(q^{n+1}-q^j) \sum\limits_{\ell=0}^{j-1} (-1)^\ell r^{ j-2\ell-1 +(j-\ell-1)^2+\ell^2} \gauss{i}{\ell}\gauss{n-i}{j-\ell-1}\\
&+\bigl((q-1)q^{j}+q^i - q^{n-i+1}\bigr) \sum\limits_{\ell=0}^{j} (-1)^\ell r^{ j-2\ell +(j-\ell)^2+\ell^2} \gauss{i}{\ell}\gauss{n-i}{j-\ell}\\
&+(q^{j+1}-1) \sum\limits_{\ell=0}^{j+1} (-1)^\ell r^{ j-2\ell+1 + (j-\ell+1)^2+\ell^2} \gauss{i}{\ell}\gauss{n-i}{j-\ell+1}.
\end{align*}
Now it suffices to show that the generating function $Z(t)=\sum_{j=0}^\infty z_j t^j$ vanishes. We first express $Z(t)$ in terms of the
polynomials $E_m(t)$ defined in Section~\ref{gaussian}.
Using Proposition~\ref{prop2.2}(iii),
we are able to rewrite our generating function as $Z(t)=\Sigma_1+\Sigma_2+\cdots+\Sigma_6$ where
\begingroup
\allowdisplaybreaks
\begin{align*}
\Sigma_1&=q^{n+1}\sum\limits_{j=0}^{\infty}\sum\limits_{\ell=0}^{j-1}(-1)^\ell r^{ j-2\ell-1 +(j-\ell-1)^2+\ell^2}\gauss{i}{\ell}\gauss{n-i}{j-\ell-1}t^j\\
&=q^{n+1}tE_i(-t)E_{n-i}(qt);\\
\Sigma_2&=-\sum\limits_{j=0}^{\infty}q^j\sum\limits_{\ell=0}^{j-1}(-1)^\ell r^{ j-2\ell-1 +(j-\ell-1)^2+\ell^2}\gauss{i}{\ell}\gauss{n-i}{j-\ell-1}t^j\\
&=-qtE_i(-qt)E_{n-i}(q^2t);\\
\Sigma_3&=(q-1)\sum\limits_{j=0}^{\infty}q^j\sum\limits_{\ell=0}^{j+1}(-1)^\ell r^{ j-2\ell +(j-\ell)^2+\ell^2}\gauss{i}{\ell}\gauss{n-i}{j-\ell}t^j\\
&=(q-1)E_i(-qt)E_{n-i}(q^2t);\\
\Sigma_4&=(q^i - q^{n-i+1})\sum\limits_{j=0}^{\infty}\sum\limits_{\ell=0}^{j+1} (-1)^\ell r^{ j-2\ell+(j-\ell)^2+\ell^2}\gauss{i}{\ell}\gauss{n-i}{j-\ell}t^j\\
&=(q^i-q^{n-i+1})E_i(-t)E_{n-i}(qt);\\
\Sigma_5&=\sum\limits_{j=0}^{\infty}q^{j+1}\sum\limits_{\ell=0}^j(-1)^\ell r^{ j-2\ell+1 +(j-\ell+1)^2+\ell^2}\gauss{i}{\ell}\gauss{n-i}{j-\ell+1}t^j\\
&=\tfrac{1}{t}E_i(-qt)E_{n-i}(q^2t);\\
\Sigma_6&=-\sum\limits_{j=0}^{\infty}\sum\limits_{\ell=0}^j (-1)^\ell r^{ j-2\ell+1 +(j-\ell+1)^2+\ell^2}\gauss{i}{\ell}\gauss{n-i}{j-\ell+1}t^j\\
&=-\tfrac{1}{t}E_i(-t)E_{n-i}(qt).
\end{align*}
\endgroup
Using Proposition~\ref{prop2.2}(i,ii), we find 
\begin{align*}
Z(t)&{}=\Sigma_1+\Sigma_2+ \Sigma_3+\Sigma_4+\Sigma_5+\Sigma_6\\
&{}=\Bigl(q^{n+1}t+q^i-q^{n-i+1}-\tfrac{1}{t}\\
&\qquad{}+\tfrac{(1-q^it)}{(1-t)}\tfrac{(1+q^{n-i+1}t)}{(1+qt)}\bigl(-qt+q-1+\tfrac{1}{t}\bigr)\Bigr) E_i(-t) E_{n-i}(qt)=0
\end{align*}
as required.


The strategy for showing $\Ptwo\Mtwo=\Dtwo\Ptwo$ is very similar but the details are sufficiently different that we provide the details here.
Fix $i$ and let $v_i=(\Ptwo_{i0},\Ptwo_{i1},\ldots,\Ptwo_{in})$.
We must show that $v_i\Mtwo=\Ptwo_{i1}v_i$. In particular, we need to show the following recurrence holds for all $j$:
\[b_{j-1}\sum\limits_{\ell=0}^{j-1} (-1)^\ell r^{ (j-1-\ell)^2+\ell^2} \gauss{i}{\ell}\gauss{n-i}{j-1-\ell} + c_{j+1}\sum\limits_{\ell=0}^{j+1} (-1)^\ell r^{ (j+1-\ell)^2+\ell^2}  \gauss{i}{\ell}\gauss{n-i}{j+1-\ell}\]
\[=\Ptwo_{i1} \sum\limits_{\ell=0}^j (-1)^\ell r^{ (j-\ell)^2+\ell^2} \gauss{i}{\ell}\gauss{n-i}{j-\ell}.\]
Multiplying both sides by $q-1$ and substituting for $b_{j-1}, c_{j+1}$, we find this is equivalent to showing that the following is zero for all $j$: 
\begin{align*}
z_j={}&(q^{n+1}-q^j) \sum\limits_{\ell=0}^{j-1} (-1)^\ell r^{(j-\ell-1)^2+\ell^2} \gauss{i}{\ell}\gauss{n-i}{j-\ell-1}\\
&+(q^{j+1}-1) \sum\limits_{\ell=0}^{j+1} (-1)^\ell r^{(j-\ell+1)^2+\ell^2}\gauss{i}{\ell}\gauss{n-i}{j-\ell+1}\\
&-q^i(q^{2n-i}-1) \sum\limits_{\ell=0}^j (-1)^\ell r^{(j-\ell)^2+\ell^2}\gauss{i}{\ell}\gauss{n-i}{j-\ell}.
\end{align*}
Again, it suffices to show that the generating function $Z(t)=\sum_{j=0}^\infty z_j t^j$ vanishes. As before, we first
rewrite our generating function as $Z(t)=\Sigma_1+\Sigma_2+\cdots+\Sigma_6$ where
\begingroup
\allowdisplaybreaks
\begin{align*}
\Sigma_1&=q^{n+1} \sum\limits_{\ell=0}^{j-1} (-1)^\ell r^{ (j-\ell-1)^2+\ell^2} t^\ell \gauss{i}{\ell}\gauss{n-i}{j-\ell-1} t^m=q^{n+1}tE_i(-rt) E_{n-i}(rt);\\
\Sigma_2&=-q^{j} \sum\limits_{\ell=0}^{j-1} (-1)^\ell r^{ (j-\ell-1)^2+\ell^2}t^\ell \gauss{i}{\ell}\gauss{n-i}{j-\ell-1}t^m=-qtE_i(-r^3t) E_{n-i}(r^3t);\\
\Sigma_3&=q^{j+1} \sum\limits_{\ell=0}^{j+1} (-1)^\ell r^{ (j-\ell+1)^2+\ell^2}t^\ell  \gauss{i}{\ell}\gauss{n-i}{j-\ell+1}t^m=\tfrac{1}{t}E_i(-r^3t) E_{n-i}(r^3t);\\
\Sigma_4&=-\sum\limits_{\ell=0}^{j+1} (-1)^\ell r^{ (j-\ell+1)^2+\ell^2}t^\ell  \gauss{i}{\ell}\gauss{n-i}{j-\ell+1}t^m=-\tfrac{1}{t}E_i(-rt) E_{n-i}(rt);\\
\Sigma_5&=-rq^{2n} \sum\limits_{\ell=0}^j (-1)^\ell r^{ (j-\ell)^2+\ell^2}t^\ell \gauss{i}{\ell}\gauss{n-i}{j-\ell}t^m=-r q^{n-i}E_i(-rt) E_{n-i}(rt);\\
\Sigma_6&=rq^i \sum\limits_{\ell=0}^j (-1)^\ell r^{ (j-\ell)^2+\ell^2}t^\ell \gauss{i}{\ell}\gauss{n-i}{j-\ell}t^m=r q^{i}E_i(-rt) E_{n-i}(rt)
\end{align*}
\endgroup
in terms of the polynomials $E_m(t)$ defined
in Section~\ref{gaussian}. Using Proposition~\ref{prop2.2}(iii), we find 
\begin{align*}
Z(t)&{}=\Sigma_1+\Sigma_2+\Sigma_3+\Sigma_4+\Sigma_5+\Sigma_6\\
&=\Bigl(q^{n+1}t-\tfrac{1}{t}-rq^{n-i}+rq^i\\
&\qquad+\bigl(\tfrac{1}{t}-qt\bigr)\tfrac{(1-rq^it)(1+rq^{n-i}t)}{(1-rt)(1+rt)}\Bigr)E_i(-rt) E_{n-i}(rt)=0
\end{align*}
as required.\qed
\end{pf} 

\bigskip

\begin{corollary}\label{corollary8.2}
The $P$-matrix of the $Q$-polynomial scheme $\mathcal{S}$ is given by
\[\begin{cases}
P_{i,j}=P_{i,2n+1-j}=\Pone_{\frac{i}2,j}&\mbox{for $i$ even, $0\leqslant j\leqslant n$;}\\
\noalign{\smallskip}
P_{i,j}=-P_{i,2n+1-j}=\Ptwo_{\lfloor\frac{i}2\rfloor,j}&\mbox{for $i$ odd, $0\leqslant j\leqslant n$.}
\end{cases}\]
\end{corollary}

\begin{pf}
If $(v_0,\ldots,v_n)$ is a left eigenvector of $\Mone$ or $\Mtwo$, it is easily seen that
either $(v_0,\ldots,v_n, v_n,\dots,v_0)$ or  
$(v_0,\ldots, v_n,-v_n,\ldots,-v_0)$ is a left eigenvector of $L_1$, respectively.  The fact that this ordering of the eigenvalues of $L_1$ is a $Q$-polynomial ordering follows from Theorem~\ref{theorem7.2}.\qed
\end{pf}

\section{A hypothetical subscheme}
\label{subscheme}

We ask whether $\mathcal{S}$ is the extended $Q$-bipartite double (in the sense of \cite{MMW}) of a primitive $Q$-polynomial scheme.  
We investigated these parameters up to $n=20$ and found they satisfied the Krein conditions, had integral eigenvalue multiplicities and nonnegative integral $p^k_{ij}$, and satisfy the handshaking lemma for all square $q$.
This appears to give an infinite family of feasible parameters for primitive $Q$-polynomial schemes with an unbounded number of classes. 
Detailed parameters and a proof of feasibility will be given in a forthcoming paper of Eiichi Bannai and Jianmin Ma.

We give the smallest case below for which existence is unknown:
\bigskip

\[P=\begin{pmatrix} 
 1 &  \frac{r^4 + r^3 + r^2 + r}{2} &  \frac{r^4 - r^3 + r^2 - r}{2} &  \frac{r^6 + r^4}{2} &  \frac{r^6 - r^4}{2} \\
  1 &  \frac{r^3 + r^2 + r - 1}{2} &  \frac{-r^3 + r^2 - r - 1}{2} &  \frac{r^4 - r^2}{2} &  \frac{-r^4 - r^2}{2} \\
  1 &  \frac{r^2 - 1}{2} &  \frac{r^2 - 1}{2} &  -r^2 &  0 \\
  1 &  \frac{-r^2 - 1}{2} &  \frac{-r^2 - 1}{2} &  0 &  r^2 \\
  1 &  \frac{-r^3 - r^2 - r - 1}{2} &  \frac{r^3 - r^2 + r - 1}{2} &  \frac{r^4 + r^2}{2} &  \frac{-r^4 + r^2}{2}
 \end{pmatrix}\]

\[Q= \begin{pmatrix} 
 1 &  \frac{r^4 - 1}{2} &  \frac{r^6 + r^4 + r^2 + 1}{2} &  \frac{r^6 - r^4 + r^2 - 1}{2} &  \frac{r^4 + 1}{2} \\
  1 &  \frac{r^4 - 2r + 1}{2r} &  \frac{r^5 - r^4 + r - 1}{2r} &  \frac{-r^5 + r^4 - r + 1}{2r} &  \frac{-r^4 - 1}{2r} \\
  1 &  \frac{-r^4 - 2r - 1}{2r} &  \frac{r^5 + r^4 + r + 1}{2r} &  \frac{-r^5 - r^4 - r - 1}{2r} &  \frac{r^4 + 1}{2r} \\
  1 &  \frac{r^4 - 2r^2 + 1}{2r^2} &  \frac{-r^4 - 1}{r^2} &  0 &  \frac{r^4 + 1}{2r^2} \\
  1 &  \frac{-r^4 - 2r^2 - 1}{2r^2} &  0 &  \frac{r^4 + 1}{r^2} &  \frac{-r^4 - 1}{2r^2}
\end{pmatrix}\]

\[L_0= \begin{pmatrix} 
 1 &  0 &  0 &  0 &  0 \\
  0 &  1 &  0 &  0 &  0 \\
  0 &  0 &  1 &  0 &  0 \\
  0 &  0 &  0 &  1 &  0 \\
  0 &  0 &  0 &  0 &  1
\end{pmatrix}\]

\[L_1= \begin{pmatrix} 
 0 &  \frac{r^4 + r^3 + r^2 + r}{2} &  0 &  0 &  0 \\
  1 &  \frac{r^2 + 2r - 3}{4} &  \frac{r^2 - 1}{4} &  \frac{r^4 + r^3}{2} &  0 \\
  0 &  \frac{r^2 + 2r + 1}{4} &  \frac{r^2 - 1}{4} &  0 &  \frac{r^4 + r^3}{2} \\
  0 &  \frac{r^2 + 2r + 1}{2} &  0 &  \frac{r^4 + r^3 + r^2 - r - 2}{4} &  \frac{r^4 + r^3 - r^2 - r}{4} \\
  0 &  0 &  \frac{r^2 + 1}{2} &  \frac{r^4 + r^3 + r^2 + r}{4} &  \frac{r^4 + r^3 - r^2 + r - 2}{4}
\end{pmatrix}\]

\[L_2= \begin{pmatrix} 
 0 &  0 &  \frac{r^4 - r^3 + r^2 - r}{2} &  0 &  0 \\
  0 &  \frac{r^2 - 1}{4} &  \frac{r^2 - 2r + 1}{4} &  0 &  \frac{r^4 - r^3}{2} \\
  1 &  \frac{r^2 - 1}{4} &  \frac{r^2 - 2r - 3}{4} &  \frac{r^4 - r^3}{2} &  0 \\
  0 &  0 &  \frac{r^2 - 2r + 1}{2} &  \frac{r^4 - r^3 + r^2 + r - 2}{4} &  \frac{r^4 - r^3 - r^2 + r}{4} \\
  0 &  \frac{r^2 + 1}{2} &  0 &  \frac{r^4 - r^3 + r^2 - r}{4} &  \frac{r^4 - r^3 - r^2 - r - 2}{4}
\end{pmatrix}\]

\[L_3= \begin{pmatrix} 
 0 &  0 &  0 &  \frac{r^6 + r^4}{2} &  0 \\
  0 &  \frac{r^4 + r^3}{2} &  0 &  \frac{r^6 + r^4 - 2r^3}{4} &  \frac{r^6 - r^4}{4} \\
  0 &  0 &  \frac{r^4 - r^3}{2} &  \frac{r^6 + r^4 + 2r^3}{4} &  \frac{r^6 - r^4}{4} \\
  1 &  \frac{r^4 + r^3 + r^2 - r - 2}{4} &  \frac{r^4 - r^3 + r^2 + r - 2}{4} &  \frac{r^6 + 2r^4 - 3r^2}{4} &  \frac{r^6 - 2r^4 + r^2}{4} \\
  0 &  \frac{r^4 + r^3 + r^2 + r}{4} &  \frac{r^4 - r^3 + r^2 - r}{4} &  \frac{r^6 - r^2}{4} &  \frac{r^6 - r^2}{4}
\end{pmatrix}\]

\[L_4= \begin{pmatrix} 
 0 &  0 &  0 &  0 &  \frac{r^6 - r^4}{2} \\
  0 &  0 &  \frac{r^4 - r^3}{2} &  \frac{r^6 - r^4}{4} &  \frac{r^6 - 3r^4 + 2r^3}{4} \\
  0 &  \frac{r^4 + r^3}{2} &  0 &  \frac{r^6 - r^4}{4} &  \frac{r^6 - 3r^4 - 2r^3}{4} \\
  0 &  \frac{r^4 + r^3 - r^2 - r}{4} &  \frac{r^4 - r^3 - r^2 + r}{4} &  \frac{r^6 - 2r^4 + r^2}{4} &  \frac{r^6 - 2r^4 + r^2}{4} \\
  1 &  \frac{r^4 + r^3 - r^2 + r - 2}{4} &  \frac{r^4 - r^3 - r^2 - r - 2}{4} &  \frac{r^6 - r^2}{4} &  \frac{r^6 - 4r^4 + 3r^2}{4}
\end{pmatrix}\]

\[L_0^*=\begin{pmatrix}
1&0&0&0&0\\
0&1&0&0&0\\
0&0&1&0&0\\
0&0&0&1&0\\
0&0&0&0&1
\end{pmatrix}\]

\[L_1^*=\begin{pmatrix}
0&\frac{r^4 - 1}{2}&0&0&0\\
1&\frac{r^6 - 5r^4 - 3r^2 - 1}{2(r^4 + r^2)} &\frac{r^8 +2r^4 + 1}{2(r^4 + r^2)} &0&0\\
0&\frac{r^6 - r^4 + r^2 - 1}{2(r^4 + r^2)} &\frac{r^8 - 4r^2 + 3}{4(r^4 + r^2)} &\frac{r^6 - r^4 + r^2 - 1}{4r^2}&0\\
0&0&\frac{r^6 + r^4 + r^2 + 1}{4r^2}&\frac{r^6 - 3r^4 - 3r^2 - 3}{4r^2}&\frac{r^4 + 1}{2r^2}\\
0&0&0&\frac{r^6 - r^4 + r^2 - 1}{2r^2}&\frac{r^4 - 2r^2 + 1}{2r^2}
\end{pmatrix}\]

\[L_2^*=\begin{pmatrix}
0&0&\frac{r^6 + r^4 + r^2 + 1}{2}&0&0\\
0&\frac{r^8 + 2r^4 + 1}{2(r^4 + r^2)} &\frac{r^{10} + r^8 + 2r^6 - 2r^4 + r^2 - 3}{4(r^4 + r^2)} &\frac{r^8 + 2r^4 + 1}{4r^2}&0\\
1&\frac{r^8 - 4r^2 + 3}{4(r^4 + r^2)} &\frac{r^{10} + 3r^8 + 2r^6 - 2r^4 + r^2 - 5}{4(r^4 + r^2)} &\frac{r^8 - 2r^6 + 2r^4 - 2r^2 + 1}{4r^2}&
\frac{r^6 + r^4 + r^2 + 1}{4r^2}\\
0&\frac{r^6 + r^4 + r^2 + 1}{4r^2}&\frac{r^8 - 1}{4r^2}&\frac{r^8 + 2r^4 + 1}{4r^2}&\frac{r^6 - r^4 + r^2 - 1}{4r^2}\\
0&0&\frac{r^8 + 2r^6 + 2r^4 + 2r^2 + 1}{4r^2}&\frac{r^8 - 2r^6 + 2r^4 - 2r^2 + 1}{4r^2}&\frac{r^6 - r^4 + r^2 - 1}{2r^2}
\end{pmatrix}\]

\[L_3^*=\begin{pmatrix}
0&0&0&\frac{r^6 - r^4 + r^2 - 1}{2}&0\\
0&0&\frac{r^8 + 2r^4 + 1}{4r^2}&\frac{r^{10} - 3r^8 - 2r^6 - 6r^4 - 3r^2 - 3}{4(r^4 + r^2)} &\frac{r^8 + 2r^4 + 1}{2(r^4 + r^2)} \\
0&\frac{r^6 - r^4 + r^2 - 1}{4r^2}&\frac{r^8 - 2r^6 + 2r^4 - 2r^2 + 1}{4r^2}&\frac{r^{10} - r^8 + 2r^6 - 2r^4 + r^2 - 1}{4(r^4 + r^2)} &\frac{r^8 - 2r^6 + 2r^4 - 2r^2 + 1}{4(r^4 + r^2)} \\
1&\frac{r^6 - 3r^4 - 3r^2 - 3}{4r^2}&\frac{r^8 + 2r^4 + 1}{4r^2}&\frac{r^8 - 4r^6 + 4r^4 - 4r^2 + 3}{4r^2}&\frac{r^6 - r^4 + r^2 - 1}{4r^2}\\
0&\frac{r^6 - r^4 + r^2 - 1}{2r^2}&\frac{r^8 - 2r^6 + 2r^4 - 2r^2 + 1}{4r^2}&\frac{r^8 - 2r^6 + 2r^4 - 2r^2 + 1}{4r^2}&0
\end{pmatrix}\]

\[L_4^*=\begin{pmatrix}
0&0&0&0&\frac{r^4+1}{2}\\
0&0&0&\frac{r^8 + 2r^4 + 1}{2(r^4 + r^2)} &\frac{r^6 - r^4 + r^2 - 1}{2(r^4 + r^2)} \\
0&0&\frac{r^6 + r^4 + r^2 + 1}{4r^2}&\frac{r^8 - 2r^6 + 2r^4 - 2r^2 + 1}{4(r^4 + r^2)} &\frac{r^6 - r^4 + r^2 - 1}{2(r^4 + r^2)} \\
0&\frac{r^4 + 1}{2r^2}&\frac{r^6 - r^4 + r^2 - 1}{4r^2}&\frac{r^6 - r^4 + r^2 - 1}{4r^2}&0\\
1&\frac{r^4 - 2r^2 + 1}{2r^2}&\frac{r^6 - r^4 + r^2 - 1}{2r^2}&0&0
\end{pmatrix}\]

\section*{Acknowledgements}

This research was supported in part by NSF grant DMS-1400281. The software package
MAGMA~\cite{Magma} was instrumental in finding the initial examples 
of these schemes, as well as in computations related to verifying our proofs.

\end{document}